\documentclass[reprint,amsmath,amssymb,aps,prb,superscriptaddress,]{revtex4-2}

\usepackage{fancyhdr}
\usepackage[english]{babel}

\usepackage[T1]{fontenc}
\usepackage[dvipsnames]{xcolor}
\usepackage[compact]{titlesec}
\usepackage{mathtools}
\usepackage{amsthm}
\usepackage{amsfonts}
\usepackage{mathtools}
\usepackage{url}
\usepackage{float}
\usepackage{bm}
\usepackage{subfiles} 
\usepackage{multirow}
\usepackage[ruled,vlined,algosection,linesnumbered]{algorithm2e}
\usepackage{hyperref}
\definecolor{GPT}{RGB}{10,50,201}
\definecolor{HumanNeedtoReRead}{RGB}{250,15,200}
\definecolor{cream}{RGB}{222,217,201}
\newcommand{\angstrom}{\mbox{\normalfont\AA}}

\newcommand{\distp}[0]{\func{p}}

\newcommand{\fat}[1]{\ifmmode\bm{#1}\else\textbf{#1}\fi}

\newcommand{\set}[1]{\mathbb{#1}}

\newcommand{\vect}[1]{\fat{#1}}

\newcommand{\tens}[1]{\mathcal{#1}}

\newcommand{\func}[1]{\textsf{#1}}

\newcommand{\vx}[0]{\vect{x}}            

\def\todecimal(#1,#2,#3,#4){\text{b2f}(#1,\,#2:[#3,\,#4])}
\def\todecimal(#1,#2,#3,#4){\text{b2f}(\{#1,\,#2\}, [#3,\,#4])}

\def\textidx#1_#2{\ensuremath{{#1}_{\text{#2}}}}
\def\ie{\textit{i.\,e.}\xspace}
\def\eg{\textit{e.\,g.}\xspace}

\def\vlist(#1,#2,#3){\{ #1,\,  #2,\, #3 \}}
\def\setn[#1]{\ensuremath{[#1]}}
\def\setn[#1]{\ensuremath{\lfloor#1\rceil}}
\def\eqdef{:=}

\def\rem#1{}
\DeclareMathOperator*{\argmin}{arg\,min}
\DeclareMathOperator*{\abs}{abs}

\begin{document}

\preprint{123-QED}

\title{Global Optimization of Atomic Clusters via Physically-Constrained Tensor Train Decomposition}
\author{Konstantin Sozykin}
\email{[konstantin.sozykin, ko.sozykin]@skoltech.ru} \email{mail@ksozykin.ru}
\affiliation{Skolkovo Institute of Science and Technology, Bolshoy Boulevard 30, Moscow, 143026, Russia}

\author{Nikita Rybin}
\affiliation{Skolkovo Institute of Science and Technology, Bolshoy Boulevard 30, Moscow, 143026, Russia}
\affiliation{Digital Materials LLC, Kutuzovskaya Street 4A, Odintsovo, Moscow Region, 143001, Russia}

\author{Andrei Chertkov}
\affiliation{Skolkovo Institute of Science and Technology, Bolshoy Boulevard 30, Moscow, 143026, Russia}
\affiliation{Artificial Intelligence Research Institute (AIRI), Presnenskaya Naberezhnaya 6, bld. 2, Moscow, Russia}

\author{Anh-Huy Phan}
\affiliation{Skolkovo Institute of Science and Technology, Bolshoy Boulevard 30, Moscow, 143026, Russia}

\author{Ivan Oseledets}
\affiliation{Skolkovo Institute of Science and Technology, Bolshoy Boulevard 30, Moscow, 143026, Russia}
\affiliation{Artificial Intelligence Research Institute (AIRI), Presnenskaya Naberezhnaya 6, bld. 2, Moscow, Russia}

\author{Alexander Shapeev}
\affiliation{Skolkovo Institute of Science and Technology, Bolshoy Boulevard 30, Moscow, 143026, Russia}
\affiliation{Digital Materials LLC, Kutuzovskaya Street 4A, Odintsovo, Moscow Region, 143001, Russia}

\author{Ivan Novikov}
\affiliation{HSE University, Faculty of Computer Science, Pokrovsky Boulevard 11, Moscow, 109028, Russia}

\author{Gleb Ryzhakov}
\affiliation{Central University, Ulitsa Gasheka 7, Moscow, 125212, Russia}

\begin{abstract}
    The global optimization of atomic clusters represents a fundamental challenge in computational chemistry and materials science due to the exponential growth of local minima with system size (i.e., the curse of dimensionality). 
We introduce a novel framework that overcomes this limitation by exploiting the low-rank structure of potential energy surfaces through Tensor Train (TT) decomposition.
Our approach combines two complementary TT-based strategies: the algebraic TTOpt method, which utilizes maximum volume sampling, and the probabilistic PROTES method, which employs generative sampling.
A key innovation is the development of physically-constrained encoding schemes that incorporate molecular constraints directly into the discretization process.
We demonstrate the efficacy of our method by identifying global minima of Lennard-Jones clusters containing up to 45 atoms. Furthermore, we establish its practical applicability to real-world systems by optimizing 20-atom carbon clusters using a machine-learned Moment Tensor Potential, achieving geometries consistent with quantum-accurate simulations.
This work establishes TT-decomposition as a powerful tool for molecular structure prediction and provides a general framework adaptable to a wide range of high-dimensional optimization problems in computational material science.

\end{abstract}

\maketitle

\section{INTRODUCTION}
    \label{intro}
    Finding the ground state of an atomic cluster is a challenging problem in physics and chemistry because it requires global optimization of the cluster’s potential energy landscape.
This landscape is typically a nonconvex, high-dimensional function with many local minima, and the number of these minima grows exponentially with the cluster size.
This exponential growth, known as the ``curse of dimensionality'', makes exhaustive search methods impractical for locating the ground state~\cite{Bellman:DynamicProgramming}.

During the last decade, Tensor Train (TT) decomposition~\cite{oseledets2011tensor} has become a powerful tool for handling high-dimensional data and computations.
This method represents multidimensional arrays (tensors) as a series of smaller, interconnected tensors, which greatly reduces computational and memory requirements.
The TT-decomposition has found applications in solving differential equations~\cite{chertkov2016robust, oseledets2016black, chertkov2021solution, chertkov2023black, Matveev2024, Dolgov+2019+23+38}, approximating multidimensional integrals and parameter-dependent integrals~\cite{ballani2012tensor, litsarev2015fast, vysotsky2021tensor}, calculating multidimensional convolutions~\cite{khoromskij2010fast, rakhuba2015fast, jin2020ctnn}, approximating the Green's function of multidimensional differential equations~\cite{fernandez2022learning, erpenbeck2023tensor, shinaoka2023multiscale}, solving computational problems in the field of financial mathematics~\cite{glau2020low, richter2021solving, bayer2023pricing}, processing audio and video~\cite{ahmadiasl2021cross, lee2021qttnet, qiu2022efficient}, accelerating and compressing artificial neural networks~\cite{tjandra2017compressing, novikov2020tensor, wang2021nonlinear, ECCV2020}, and even building new machine learning algorithms directly based on the TT-decomposition~\cite{chen2019support, kour2023efficient, liu2023tensor}.
In the last couple of years, the TT-decomposition has been successfully applied to several gradient-free optimization problems~\cite{TTopt_NEURIPS2022, NEURIPS2023_02895786, sozykin2025high, chertkov2023tensorExtrema}, enabling efficient handling of high-dimensional functions in fields such as reinforcement learning~\cite{shetty2024generalized}, robotics~\cite{shetty2024tt_robotics, xue2025unifyingrobotoptimizationmonte}, brain studies~\cite{lee2024tt, pospelov2025fast}, epidemiological modeling~\cite{dolgov2025tensorcross}, and adversarial training of neural networks~\cite{chertkov2023tensor, chertkov2024translate}.

Inspired by the above works, we develop a novel framework to address the challenging task of optimizing the geometry of atomic clusters. 
In this framework, tensors are employed as a low-rank model to search for candidates of potential energy minima.
To achieve this, we utilize the generalized maximum volume principle~\cite{Goreinov2010, mikhalev2018rectangular} and probabilistic sampling techniques~\cite{dolgov2020approximation, NEURIPS2023_02895786} followed by single-shot relaxation via the L-BFGS-B algorithm~\cite{Nocedal2006Numerical}.

We primarily adopt the TT-format due to its computational efficiency and convenience.
Nevertheless, the proposed framework can be extended to other tensor decomposition schemes. 
Examples include the Canonical Polyadic Decomposition (CPD-PARAFAC)~\cite{CPD, harshman1970foundations}, Tucker decomposition~\cite{Tucker_1966}, Hierarchical Tucker decomposition~\cite{buczyska2015hackbusch, ballani2013black, ryzhakov2024black}, and Tensor Ring decomposition~\cite{TensorRing, TensorRing2}.
Additionally, tensor networks traditionally associated with quantum physics, such as Projected Entangled Pair States (PEPS)~\cite{PhysRevA.70.060302, verstraete2004} and the Multiscale Entanglement Renormalization Ansatz (MERA)~\cite{PhysRevLett.101.110501, PhysRevB.79.144108, PhysRevB.89.235113, Batselier2021}, can also be integrated into our framework.
Other relevant methods are explored in~\cite{cichocki2016tensor, cichocki2017tensor, Orus2019}.

To validate the proposed approach, which we refer to as TT-based optimization, we applied it to Lennard-Jones clusters and successfully identified global minima for systems comprising 5 to 45 atoms.
Furthermore, we evaluated the framework's performance in optimizing 20-atom carbon clusters using the learnable Moment Tensor Potential (MTP)~\cite{shapeev2016_mtp} as the interatomic interaction model.
MTP is a state-of-the-art machine learning interatomic potential that bridges the gap between computationally expensive ab initio methods and less accurate empirical potentials. 
In our experiments, the optimized structure was determined with high precision, achieving accuracy within the fitting error of the MTP model.
These results demonstrate the efficacy of TT-based optimization in solving complex energy minimization problems across various atomic systems.
\section{METHODOLOGY}
    \label{sec:methods}
    In this section, we present our methodology, including the formulation of the considered optimization problem and a general description of the TT-decomposition; a discussion of the approaches we use for discretization of the problem; details of the proposed TT-based optimization method; and a list of potentials we consider for numerical experiments.
Table~\ref{tab:notation} summarizes the notation used throughout the manuscript; relevant definitions will also be provided in context where necessary.
\begin{table}[H]
    \small
    \centering
    \caption{Summary of notation used throughout the manuscript.}
    \label{tab:notation}
    \begin{tabular}{|l|p{0.3\linewidth}|l|p{0.3\linewidth}|}
    \hline
    \textbf{Symbol} & \textbf{Description} & \textbf{Symbol} & \textbf{Description} \\
    \hline
    \multicolumn{2}{|c|}{\textbf{General}} & \multicolumn{2}{c|}{\textbf{TT-Related}} \\
    \hline
    $M$ & Number of atoms & $\tens{G}_i$, $R_i$ & TT-core and rank \\
    $\vect{x}$ & Coordinates & TT-format & Tensor Train \\
    $\Omega$ & Search domain & TTOpt & Algebraic optimizer \\
    $d$ & Problem dimension & PROTES & Probabilistic optimizer \\
    $\vect{n}$ & Multi-index & & \\
    $E(\vect{x})$ & Energy function & & \\
    $\tens{E}[\vect{n}]$ & Energy tensor & & \\
    \hline
    \multicolumn{2}{|c|}{\textbf{Encodings}} & \multicolumn{2}{c|}{\textbf{Initialization}} \\
    \hline
    DE & Direct (Cartesian) & Agn & Agnostic (uniform) \\
    SR & Simple relative & PhC & Physically-constrained \\
    CR & Constrained relative & & \\
    Bit & Binary/p-ary & & \\
    \hline
    \multicolumn{4}{|c|}{\textbf{Metrics}} \\
    \hline
    PC & PROTES calls & LCT, LCL & Local calls (total, last) \\
    $SR_t$ & Success rate (\%) & RE & Relative error \\
    \hline
    \end{tabular}
\end{table}

\subsection{Optimization of the potential energy and the low-rank tensor representation}
    We denote the potential energy of a cluster by $E(\vect{x})$, where $\vect{x}$ is a set of particle (e.g., atom) positions.
Optimization algorithms attempt to find the optimal input $\vect{x}$ for the function $E(\vect{x})$ by solving:
\begin{equation}\label{eq:opt_def}
\begin{split}
\min_{ \vect{x} } & \quad E(\vect{x})
    \quad  subject \: to \\
\quad H_{k}(x) = 0 &
    \quad k=0,1,\ldots,K-1, \\
\quad G_{j}(x) \leq 0 & \quad j=0,1,\ldots,J-1, 
\end{split}
\end{equation}
over a specific feasible set $\vect{x}$ with possible constrains given usually as a set of equations $\{H_{k}(x)\}$, inequalities $\{G_{j}(x)\}$, and a specific domain 
$$
\Omega =
    (a_x,b_x) \times 
    (a_y,b_y) \times
    (a_z,b_z),
$$
in which we search for the optimal particle positions. 
The problem~\eqref{eq:opt_def} is difficult to solve in this formulation due to the exponential growth in the number of local minima as the number of particles increases~\cite{CSP_book_intro, FIRE}.

To simplify the solution of the original problem, we propose the following approach.
We represent the potential energy of a cluster, which is a function of positions $\vect{x}$ of $M$~particles, as a $d$-dimensional tensor $\tens{E}$, where $d=3M$. 
Note that the dimensionality of the discrete tensor can be lower depending on the encoding (the method by which we translate the positions of particles into their discrete counterparts).
One can also reduce the dimensionality by accounting for the symmetry of the potential energy with respect to rotations, translations, etc.

Here, we define a $d$-dimensional tensor $\tens{E}$ as a multidimensional array, i.e., a set of real numbers indexed by multi-indices $[n_1, n_2, \ldots, n_{d}]$, with each index $n_i=0,1,\ldots,N_i-1$ ($i=1,2,\ldots,d$), where $N_i$ is the number of possible particle positions.
We can think of a tensor as a real-valued function of $d$~variables, each of which takes values from a discrete domain.
For example, for the energy of an atomic structure, the $i$-th input argument corresponding to the $x$-coordinate of some atom from the original vector belongs to the following discrete set:
\begin{equation}\label{eq:basic_grid}
\Phi_x=\{
a_x,\,
a_x+h_i,\,
a_x+2h_i,\,
\ldots,\,
a_x+(N_i-1)h_i,\,
b_x
\},
\end{equation}
where we impose a grid on the interval $[a_x,\,b_x]$ with constant step size $h_i = \frac{b_x-a_x}{N_i}$.

\begin{figure}[t!]
    \centering
    \includegraphics
        [width=0.95\linewidth]
        {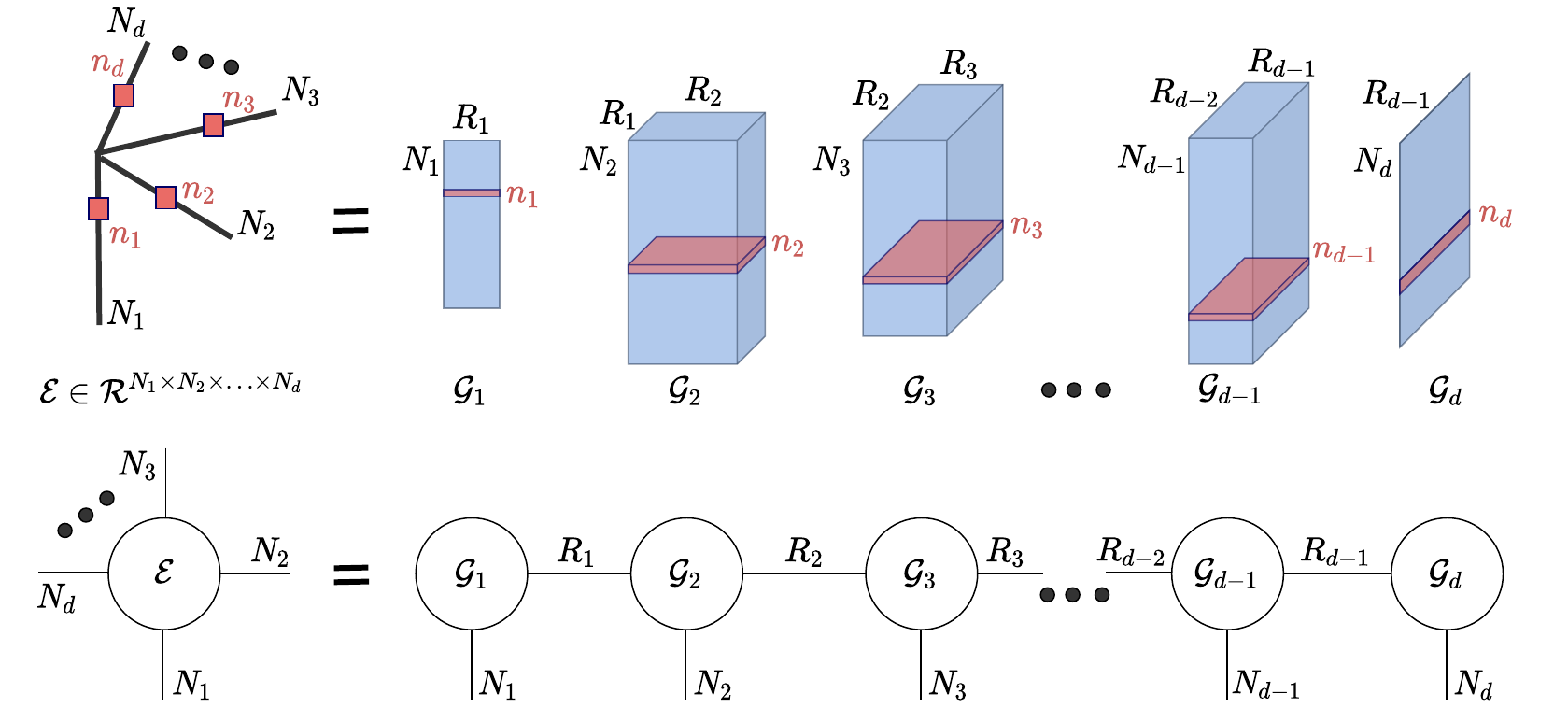}
    \caption{
        Schematic representation of the TT-decomposition.
        The procedure for calculating element $[n_1, n_2, \ldots, n_d]$ using a tensor in the TT-format is presented at the top, and the tensor diagram corresponding to the TT-decomposition is demonstrated at the bottom.
    }
    \label{fig:scheme_tt}
\end{figure}

Therefore, the original formulation~\eqref{eq:opt_def} can be efficiently redefined as a discrete optimization problem of selecting the corresponding indices of the tensor $\tens E$ that represents the discretized search space of $\vect n$:
\begin{equation}\label{eq:opt_def2}
\begin{aligned}
\min_{ \vect{n}\in\Phi} \; & \tens E[\vect n] \quad 
\text{s.t.}
\\
\quad
    \vect{n} &= [
        n_{1x},\, n_{1y},\, n_{1z},\,
        \ldots,\,
        n_{Mx},\, n_{My},\, n_{Mz}
    ],\\
    & n_{ix}\in\Phi_x,\;n_{iy}\in\Phi_y,\;n_{iz}\in\Phi_z.
\end{aligned}
\end{equation}
The advantage of such a discrete representation is the possibility of using low-rank tensor decomposition, i.e., the representation of a multidimensional array that requires far fewer elements for its storage than the total number of the given array elements.

As noted in the Introduction, we use the so-called Tensor Train (TT) as a particular low-rank tensor decomposition. 
We denote that the $d$-dimensional tensor $\tens{E} \in \mathrm{R}^{N_1 \times N_2 \times \cdots \times N_d}$ is expressed in the TT-format~\cite{oseledets2011tensor} if its elements can be represented (exactly or approximately) by the following formula (see Figure~\ref{fig:scheme_tt} for an illustration):
\begin{equation}\label{eq:tt_repr}
\begin{split}
\tens{E} [n_1, & n_2, \ldots, n_d]
= 
\sum_{r_0=1}^{R_0}
\,
\sum_{r_1=1}^{R_1}
\cdots
\sum_{r_{d}=1}^{R_{d}}
    \\
    &
    \tens{G}_1 [r_0, n_1, r_1]
    \,\,
    \tens{G}_2 [r_1, n_2, r_2]
    \,
    \times
    \\
    \times
    \,
    \cdots
    \,
    &
    \tens{G}_{d-1} [r_{d-2}, n_{d-1}, r_{d-1}]
    \,\,
    \tens{G}_d [r_{d-1}, n_d, r_d],
\end{split}
\end{equation}
where $[n_1, n_2, \ldots, n_d]$ is a multi-index (with $n_i = 0, 1, \ldots, N_i-1$ for $i = 1, 2, \ldots, d$), and the integers $R_{0}, R_{1}, \ldots, R_{d}$ (with the convention $R_{0} = R_{d} = 1$) are called TT-ranks which determine the accuracy of approximation.
The three-dimensional tensors $\tens{G}_i \in \mathbb{R}^{R_{i-1} \times N_i \times R_i}$ (for $i = 1, 2, \ldots, d$) are usually called TT-cores (or carriages of the TT-decomposition).

We note that the total storage required for the TT-cores \( \tens{G}_1, \tens{G}_2, \dotsc, \tens{G}_d \) is at most
$$
d \times \max_{1 \leq k \leq d} N_k \times \left( \max_{0 \leq k \leq d} R_k \right)^2
$$
memory cells.
Consequently, when the TT-ranks are bounded, the TT-approximation circumvents the curse of dimensionality. 
Additionally, fundamental linear algebra operations such as computing norms, differentiation, integration, solving linear system, etc., can be efficiently implemented in the TT-format with polynomial complexity relative to both the dimensionality and the mode sizes~\cite{cichocki2016tensor, cichocki2017tensor}.

\subsection{Discretizations and grid projections}\label{seq:grids}
    Since our top-level algorithm is discrete, we need a mapping of continuous particle positions into a discrete vector indices. 
We propose several such mappings, which are presented below.
We start the description with a na\"{\i}ve algorithm that simply encodes the three-dimensional Cartesian coordinates of each cluster sequentially by putting them on the predefined mesh.
We then use properties of the potential energy of the particle system, such as translational and rotational invariance, to reduce the number of the used parameters. In addition, since we are searching for a minimum energy, we do not need to have configurations that have very high energy. 
Thus, we can immediately set, for example, the minimum distance between neighboring (in the coding order) particles. This approach also allows us to reduce the dimensionality of the parameter space.

\subsubsection{Direct encoding}

Let us consider, as an example, the $x$-coordinate of the $i$-th particle (for any $i = 1, 2, \ldots, M$), which lies in the $[a_x, b_x]$ interval.
In a discrete representation, a coordinate can take one of the values from the set $\Phi_x$  defined in~\eqref{eq:basic_grid}.
At the encoding stage, we are looking for the index~$n$ of the closest point of the mesh~$x_n \in \Phi_x$ to the point~$x$:
$$
x_n=\argmin_{x_i\in\Phi_x} |x_i-x|.
$$
To find 
the index~$n$ of this point explicitly,
we define the mapping function:
\begin{equation}\label{eq:f2i}
f2i(x,\,[a,b],\,N):=\left\lfloor
\frac{N (x -a)}{b-a}+1/2
\right\rfloor,
\end{equation}
where $\lfloor q \rfloor$ denotes greatest integer less than or equal to~$q$. 
Then the translation of a continuous coordinate $x$ into a discrete index $n$ occurs according to the formula:
$$
n = f2i(x,[a_x,b_x],N_i).
$$
And vise-versa: the real value of the mesh point~$x_n$ corresponding to the index value $n$
can be obtained as
$$
x_n = i2f(n,[a_x,b_x],N_i),
$$
where we define the mapping function:
$$
    i2f(n,\,[a,b],\,N):=\frac{n}{N-1}(b-a)+a.
$$

On the next step, we assign the obtained index~$n$ with the corresponding input index component of the tensor.
Repeating this procedure for all components of all~$M$ particles, we get a set of integers, which we interpret as the input index of our tensor:
$$
\def\xn#1{\underbrace{n_{x#1},n_{y#1},n_{z#1}}_{\text{particle }#1}}
I=\{\xn1,\xn2,\ldots,\xn M\},
$$
where $n_{xi}$ is the index for $x$-coordinate of the $i$-th particle (the same is for $y$ and $z$),
obtained by the described procedure.

\rem{
At the encoding stage we translate each of the three-dimensional coordinates into a $p$-number system.
Thus, we obtain a vector~$V\in\setn[p]^{3\cdot q\cdot N}$ (we define $\setn[p]:=\{0,1,\ldots,p-1\}$) which represents the encoded positions of~$N$ atoms.}

\rem{
Let us take a closer look at the coding procedure.
A~vector~$v=\{v_1,\ldots,v_q\}$ of length~$q$, which represents one of the spatial coordinates ($x$, $y$ or $z$) of some particle, 
has elements of integer numbers from the range $v_i=0,\ldots,p-1$. 
In the case of the binary number system, its components are only zeros and ones.
Treating the elements of the given vector~$v$ as bits of the $p$-ary number system, we construct an integer~$K$ based on all elements of the vector~$v$ as follows:
\begin{equation*}
K=v_1+v_2p+v_3p^2+\cdots+v_qp^{q-1}.
\end{equation*}
The minimum value that~$K$ can take is zero (in the case where all elements of the vector~$v$ are zero),
and the maximum is~$p^q-1$ (reached when all elements of the vector~$v$ are~$(p-1)$).
We map all possible values of~$K$ so that they cover the given interval~$[a,b]$ uniformly, namely, finally the coordinate of the particle is equal to:
\begin{equation*}
x=K\frac{b-a}{p^q-1}+a.
\end{equation*}
}


\subsubsection{Relative encodings}

To enhance the procedure and adopt a more chemistry-compatible representation, we employ relative encoding -- a method conceptually aligned with the Z-matrix framework \cite{Parsons2005}. 
Firstly, each subsequent atom is encoded by its relative position in relation to one of the previous ones.
The index of this ``previous ones'' can be considered as optimized variables.
Secondly, the relative position is encoded in spherical coordinates, which allows us to set the value of the distance $r$ to the previous atom within rather narrow limits, \eg, in the interval $[1, 1.2]$.
Finally, we additionally assume that the initial atom is at the origin of the coordinates, and we place the next one after on the $z$-axis.
Any position of atoms can be reduced to such a position by translation and rotation, which does not change the interaction energy of atoms.

\paragraph{Simple relative encoding}
In this case, the resulting set of input indices will be defined as follows:
\def\ooo#1{\underbrace{%
        c_{#1},n_{r_{#1}},n_{\theta_{#1}},n_{\phi_{#1}}}_{\text{particle {#1}}%
    }}
\begin{multline}\label{ref_enc}   
I_r = \{\underbrace{n_{r_2}}_{\text{particle 2}}
    ,\ooo3
    ,\ooo4,\ldots,
    \\
    \ooo M
\}.
\end{multline}
Note, that here
$$
c_3=1,2;\quad
c_4=1,2,3;\quad
\ldots;\quad
c_M=1,2,\ldots,M-1.
$$
We use these indices to restore real coordinates of the $i$-th ($i=3,4,\ldots,M$) particle as follows:
\begin{equation*}
\begin{cases}
    x_i = x_{c_i} + r_i\cos\phi_i\sin\theta_i, \\
    y_i = y_{c_i} + r_i\sin\phi_i\sin\theta_i, \\
    z_i = z_{c_i} + r_i\cos\theta_i,
\end{cases}
\end{equation*}
where $c_i$ is the index of the previous particle, and $(r_i, \theta_i, \phi_i)$ are decoded spherical coordinates:
\begin{align}
    r_i
        &
        =i2f(n_{r_i},[r_{\text{min}},r_{\text{max}}],N_{r}),\nonumber\\
    \theta_i
        &
        =i2f(n_{\theta_i},[0,\pi], N_{\theta}),\nonumber\\
    \phi_i
        &
        = i2f(n_{\phi_i},[0,2\pi], N_{\phi}).
\end{align}
Additionally,
for the first particle ($i=1$) we have $x_1=y_1=z_1=0$;
for the second particle ($i=2$): $x_2=y_2=0$, and $z_2=i2f(n_{r_2},[r_{\text{min}},r_{\text{max}}],N_{r})$.

\emph{The advantage} of this type of coding is that we know for sure that each point is separated from another point by a specific distance that is neither too large nor too small (in the example, from~$1$ to~$1.2$), rather than being scattered infinitely. 
At the same time, the whole structure of the points themselves can be quite large in space. 
The second advantage is that the accuracy of the grid does not depend on the number of points: if we take many points, we have to expand the limits (the cells $[-L,L]$ in which the coordinates lie), \ie, the larger number of particles~$M$ is, the larger~$L$ is. 
Thus, the accuracy of the grid $h\sim L/N$ will decrease proportionally with a fixed number of mesh points~$N$.
In the presented coding, we have fixed intervals ($[1, 1.2]$ for $r$, $[0, 2\pi)$ for $\phi$,  $[0, \pi)$ for $\theta$) that do not depend on the number of points, and therefore the accuracy of the grid does not decrease with the number of points.

\paragraph{Relative encoding with constant distances}
To further simplify the optimization procedure, relative coding with \emph{constant distances} can be considered.
In this case, we assume that~$r_i$ is the same for all particles, which means that only one variable can be used in the encoding, for example, the first one:
\def\ooo#1{\underbrace{
        c_#1,n_{\theta_{#1}},n_{\phi_{#1}}
    }_{\text{particle #1}}}
$$
I_{r} = \{
    n_r,
    \ooo3,
    \ooo4,
    \ldots, \ooo M
\},
$$
where $n_r=0,1,\ldots,N_r$, and it is rational to take the number of discrete points sufficiently large, for example $N_r=2N_\theta$.
With this type of encoding, we reduce the length of the vector used to encode the structure by a factor of approximately~$4/3$.
Furthermore, it is more consistent with physical observations that each point has a neighbor at a fixed (for a given structure) distance.

\paragraph{Relative encoding with constant distances and restrictions on angles}
We can go further and prevent the distance between a given particle and its predecessor from being too small.
We achieve this by modifying the previous algorithms 
to consider angles in the spherical coordinate system as relative, while constraining the minimum value of the angle $\theta$.
Thus, we come to the relative encoding with restriction on distances.

Firstly, we decode angles (and~$r$ if we consider it different for different particles) for each particle with number $i>2$ (array index numeration starts from 1):
\begin{align}
\textidx\theta_{rel} &= i2f(n_{\theta_i},[0,\pi], N_{\theta}),\\
\textidx\phi_{rel}   &= i2f(n_{\phi_i},[0,\pi], N_{\phi}).
\end{align}
Secondly, we consider the matrix of the composition of rotations along the $z$ and $y$ axes:
\begin{equation*}
R_{ZY}(\theta,\,\phi)
=
\def\phic{\textidx\phi_{rel}}
\def\phic{\phi}
\begin{pmatrix}
\cos\phic & -\sin\phic & 0\\
\sin\phic & \cos\phic & 0\\
0 & 0 & 1\\ 
\end{pmatrix}
\times\\
\def\phic{\textidx\theta_{rel}}
\def\phic{\theta}
\begin{pmatrix}
\cos\phic & 0 & \sin\phic\\
 0 & 1 & 0\\
-\sin\phic & 0 & \cos\phic\\ 
\end{pmatrix}.
\end{equation*}

In the next step, we assign a product of such the matrices with each particle using the decoded values of angles~$\textidx\theta_{rel}$ and~$\textidx\phi_{rel}$:
$$
R_i=R_c \times R_{ZY}(\textidx\theta_{rel},\,\textidx\phi_{rel}).
$$
For the first two particles, 
which have coordinates $x[1] = \vlist(0,0,0)$ and $x[2]=\vlist(0,0,r)$, respectively, we have
\begin{equation}
R_1\eqdef 
R_{ZY}(\pi,\,0),
\quad
R_2\eqdef 
R_{ZY}(0,\,0).
\end{equation}

Finally, 
to obtain the 3D coordinates of the current particle, 
we combine these matrices for the current particle and its predecessor, and use the predecessor's coordinate as a relative position:
\begin{equation*}
    \{\hat{x}_i, \hat{y}_i, \hat{z}_i\} =
    \{x_i, y_i, z_i\} +
    r
    \cdot
    (R_{i}\times\vlist(0,0,1)^T).
    \label{eq:rel_p}
\end{equation*}
The advantage of this type of encoding is that we can halve the number of sampling points for angle 
$\theta$, thereby simplifying the search for the optimum. 
Secondly, we reduce the chances of particles ``colliding'', \ie{} the potential being close to infinity.

Summarizing the ideas presented in the encodings, we note that they not only technically reduce the number of parameters for optimization, but also set certain \emph{constrains} so that particle configurations are more physical.
Note that such constraints can be set further in a similar style, but in order not to increase the complexity of the algorithm, we limited ourselves to those described.

\subsubsection{Bit coding}
The use of low-rank tensor decomposition itself enables for solving the problem of the curse of dimensionality.
For example, for the TT-decomposition, the number of tensor elements for naive coding is $O(3M\cdot N\cdot R^2)$, where~$N$ is the average number of points on the mesh of each coordinate, and $R$ is the average TT-rank.
However, if we take a huge number of such points, \ie, if $N$ is large, the number of tensor parameters can also be large.
Here we describe a procedure for further reducing the number of parameters in the low-rank decomposition, which is based on the \emph{bit encoding} procedure of numbers~$n_i$, which are components of the tensor input.
Note that each bit coding is simply an extension of one of the coding methods described above.

Let us take a closer look at the bit encoding procedure.
We fix the base of the number system~$p>1$ and the number of bits~$q>0$ which we use to encode the particle parameter~$n_i$.
It means that the number of possible values that~$n_i$ can take is~$N_i=p^{q}$ (\ie{} $n_i=0,1,\ldots,p^{q}-1$).
We assign a~vector~$v=\{v_1,\ldots,v_q\}$ of length~$q$, 
with the number~$n_i$.
This vector
has elements of integer numbers from the range $v_i=0,1,\ldots,p-1$. 
In the case of the binary number system ($p=2$), its components are only zeros and ones.
Treating the elements of the given vector~$v$ as bits of the $p$-ary number system, we reconstruct an integer~$n_i$ based on all elements of the vector~$v$ as follows:
\begin{equation*}
n_i=v_1+v_2p+v_3p^2+\cdots+v_qp^{q-1}.
\end{equation*}
Thus, instead of a single index (input) of the tensor, which can vary from~$0$ to~$N_i-1$, 
we get $q$ indices, each of which can vary from~$0$ to~$p-1$.
The minimum value that~$n_i$ can take is zero (in the case where all elements of the vector~$v$ are zero),
and the maximum is~$p^q-1$ (reached when all elements of the vector~$v$ are~$(p-1)$).

With this encoding, we can directly get the real values~$x$ from a given interval~$[a_x,b_x]$.
To do this, we map all possible values of~$n_i$ so that they cover the given interval~$[a,b]$ uniformly, namely, the coordinate of the particle~$x$ is equal to:
\begin{equation*}
x=n_i\frac{b-a}{p^q-1}+a.
\end{equation*}
For any consistent base~$p$ and vector of bits~$v$ of the length~$p$ and real numbers $b>a>0$,
we denote the full procedure as 
\begin{multline}\label{eq:def_of_todec}
\todecimal(v,p,a,b):=\\
\bigl(v_1+v_2p+v_3p^2+\cdots+v_qp^{q-1}\bigr)\frac{b-a}{p^q-1}+a.
\end{multline}

Using this approach, we obtain for the na\"{\i}ve case
the following input vector of indices for the tensor:
\def\ub#1#2#3{\underbrace{
v_{#1},\ldots,v_{#2}
}_{#3}}
\begin{multline*}
I=\{
\ub1q{x_1},
\ub{q+1}{2q}{y_1},
\ub{2q+1}{3q}{z_1},\\
\ub{3q+1}{4q}{x_2},
\,\ldots\,,
\ub{(3M-1)q}{3Mq}{z_M}
\},
\end{multline*}
where $x_i$, $y_i$ and $z_i$ are spatial mesh coordinates of the $i$-th particle.
The density of this mesh depends on the parameters~$p$ and~$q$.
The number of low-rank parameters in the TT-decomposition is now $O(3M\cdot qp\cdot r^2)$.
Similarly, bit coding extends other encoding methods presented above; the technical details are very similar, so we will not go into them here.


\subsection{Optimization methods based on tensor trains}
    \begin{figure}[t!]
    \centering
    \includegraphics
        [width=0.95\linewidth]
        {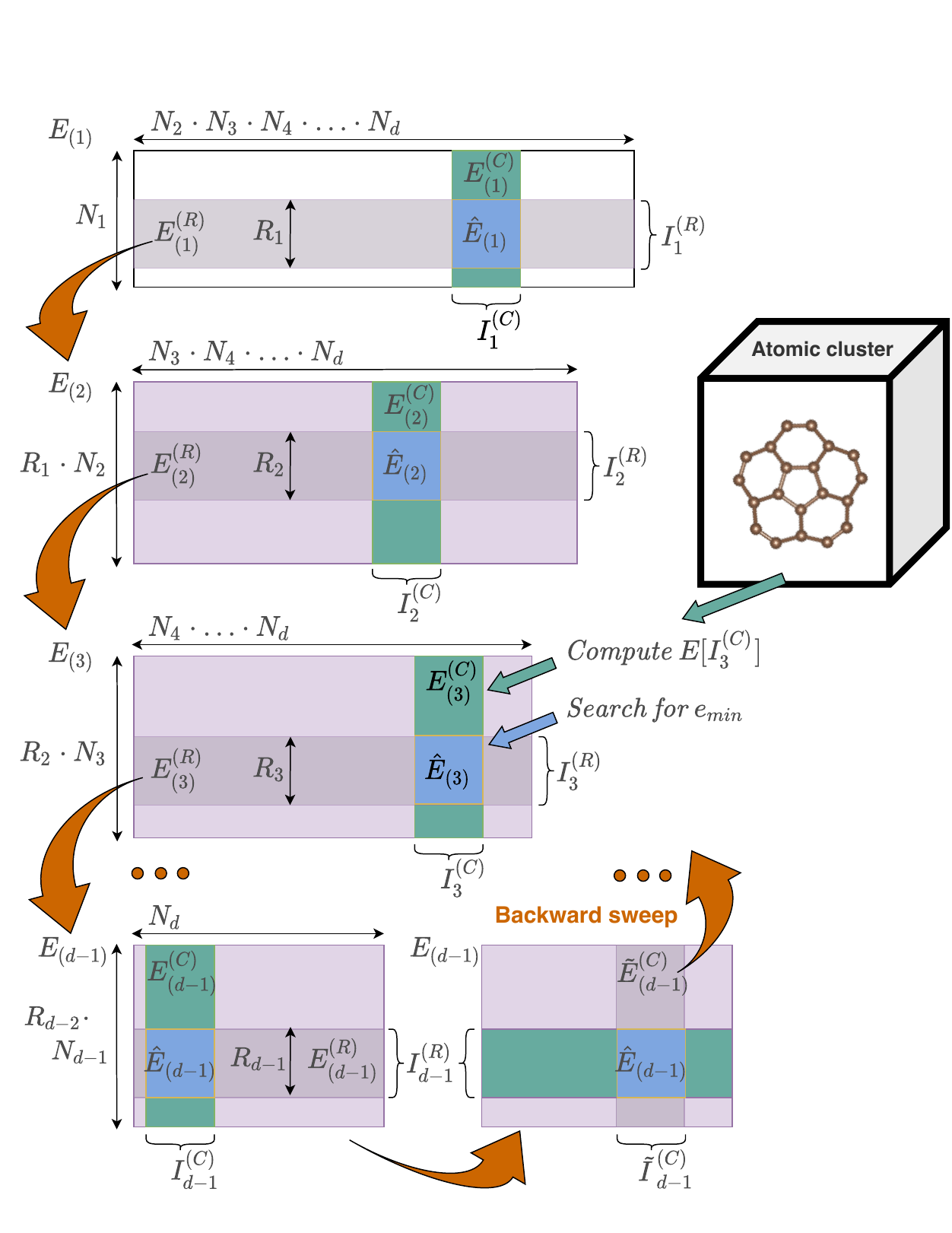}
    \caption{
        Schematic representation of the optimization method TTOpt.
    }
    \label{fig:scheme_ttopt}
\end{figure}

\begin{figure*}[t!]
    \centering
    \includegraphics
        [width=0.95\linewidth]
        {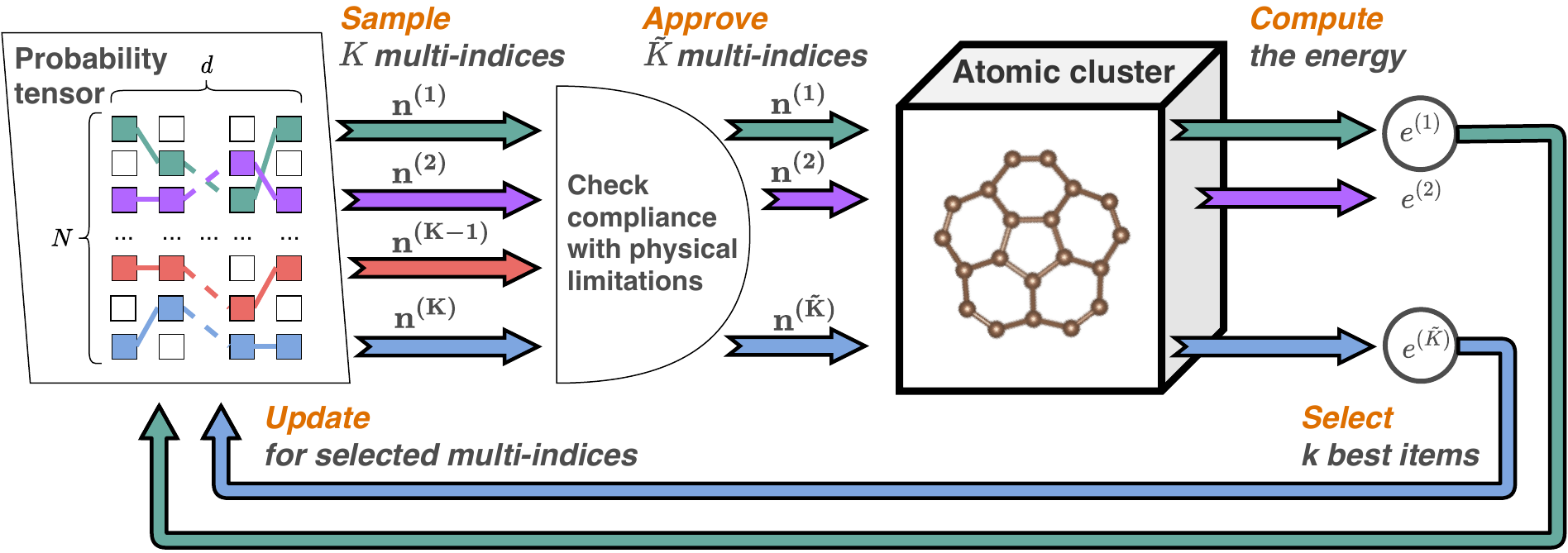}
    \caption{
        Schematic representation of the optimization method PROTES.
    }
    \label{fig:scheme_protes}
\end{figure*}

This section details two TT-based approaches adapted for high-dimensional optimization of atomic clusters: maxvol-based optimization (TTOpt)~\cite{TTopt_NEURIPS2022} and probabilistic sampling (PROTES)~\cite{NEURIPS2023_02895786}.
Both methods are applicable to the discrete minimization problem~\eqref{eq:opt_def2} for the $d$-dimensional potential energy tensor $\tens{E} \in \mathbb{R}^{N_{1} \times N_2 \times \cdots \times N_{d}}$ defined on uniform grids $\Phi_x, \Phi_y, \Phi_z$ of the form~\eqref{eq:basic_grid}.
The workflows for these methods are illustrated in Figures~\ref{fig:scheme_ttopt}, and~\ref{fig:scheme_protes}, and are discussed in detail below.

\subsubsection{Method TTOpt}

TTOpt leverages the maximum volume (``maxvol'') principle~\cite{Goreinov2010} to efficiently locate near-global minima in high-dimensional spaces.
The concept behind the maxvol algorithm is as follows.
Consider a ``tall'' matrix~$A\in\mathbb R^{N_1\times N_2}$, where the number of rows~$N_1$ is significantly larger than the number of columns~$N_2$.
The objective is to select a subset of row indices~$I$ such that the corresponding rows are as ``independent'' as possible.
This requirement is formalized by maximizing the volume of the submatrix~$B=A[I, :]$, which is defined as $\text{vol}(B)=\sqrt{\det{B^TB}}$.
The core theorem (\cite[Theorem 1]{Goreinov2010}) states that for any rank-$r$ matrix $A$, there exists an $r \times r$ submatrix~$B'$ whose element with the maximum modulus is at least $1/r^2$ of the global maximum:
$$
\max \abs B' \geq \frac1{r^2}\max \abs A.
$$
Here ``$\abs$'' denotes the element-wise absolute value function.
Under these conditions, the efficient iterative maxvol algorithm~\cite{Goreinov2010} can be employed to find such a submatrix $B'$ by maximizing its volume.
This estimate enables an exponential reduction of the search space while preserving high-quality solutions.
In practice, we assume that matrix~$A$ is of full rank ($r=N_2$), although we can select more row indices than the number of columns~$N_2$ by using the ``rectangular'' version of the maxvol algorithm~\cite{mikhalev2018rectangular}.
Thus, we expect that there is an element in such a submatrix, which is close to the maximum of the entire matrix~$A$.
Note that if we want to \emph{minimize} the black box function~$f$ rather than maximize, we do not use its values directly, but apply a strictly decreasing monotonic transformation that converts the minimum to the maximum and vice versa. In our experiments, we used the exponential transformation: $f(x)\mapsto\exp(-f(x))$.

The problem of finding the maximal (or minimal) element in a multidimensional tensor can be reduced to finding the maximal element in a matrix via tensor unfolding.
A tensor unfolding is a matrix whose rows and columns are indexed by super-indices.
For the tensor $\tens{E}$, its $l$-th ($l=1,2,\ldots, d-1$) unfolding $E_{(l)}$ is defined by its elements as
$$
    E_{(l)}[
        \overline{n_1 \dots n_l},
        \overline{n_{l+1} \dots n_d}
    ] =
    {\tens E}[n_1, \, n_2, \, \ldots, \, n_d],
$$
which are obtained by combining the tensor indices into super-indices using the row-major convention:
\begin{equation}\label{eq:def_super-indices}
\overline{n_1 \dots n_l} = \sum_{j=1}^l \left( n_j \prod_{s=1}^{j-1} N_s \right).
\end{equation}

The TTOpt algorithm (see Figure~\ref{fig:scheme_ttopt}) essentially processes the tensor unfoldings~$E_{(l)}$ sequentially, attempting to consistently identify pairs of super-indices that define a submatrix with maximum volume.
At the initial stage, we set the expected ranks of the approximation: $R_1, R_2, \ldots R_{d-1}$.
In each sweep, starting from the first core, we form a submatrix~$E^{(C)}_{(1)}$ of the unfolding~$E_{(1)}$ by selecting all rows of the unfolding and columns with indices $I^{(C)}_1$, which were obtained in previous sweeps or initialized randomly during the first cycle.
The submatrix $E^{(C)}_{(1)} = E_{(1)}[:, I^{(C)}_1] \in \set{R}^{N_1 \times R_1}$ is constructed by 
calculating all its elements using the black box function (\ie, transformed energy values).
Then, using the maxvol algorithm, we identify the maximal-volume submatrix $\hat{E}_{(1)} \in \set{R}^{R_1 \times R_1}$ within $E^{(C)}_{(1)}$ and store its row indices in the list $I^{(R)}_1$.
These indices correspond to the values of the first tensor index $n_1=1, 2, \ldots, N_1$ that provide the most information.

Subsequently, we proceed to the next core.
The indices $I^{(R)}_1$ are used to generate a column submatrix $E^{(C)}_{(2)}$ in the next unfolding matrix $E_{(2)}$.
Specifically, from each index \(i \in I^{(R)}_1\), \(N_2\) distinct indices are generated by systematically appending each permissible value \(v = 1, 2, \dots, N_2\) of the tensor's second leg to \(i\).
This yields an index set $\widetilde I^{(R)}_2$ of size $|\widetilde I^{(R)}_2|=|I^{(R)}_1|\cdot N_2=R_1\cdot N_2$.
Here, $|\cdot|$ denotes the cardinality of the set.
The submatrix $E^{(C)}_{(2)}=E_{(2)}[\widetilde I^{(R)}_2, I^{(C)}_2]\in\mathbb R^{R_1N_2\times R_2}$ is populated using the given black box.

We then apply the maxvol algorithm to the submatrix~$E^{(C)}_{(2)}$ to obtain the row index list $I^{(R)}_2\subset \widetilde I^{(R)}_2$ with length $|I^{(R)}_2|=R_2$.
We assume that the ranks are selected such that $R_1N_2\geq R_2$ holds.
Similar to the previous case, the set of indices $I^{(R)}_2$ represents the $R_2$ pairs of values for the first two tensor legs that carry the maximum amount of information.

This process is repeated until the last tensor unfolding is reached, followed by a backward pass from the last to the first unfolding.
During the reverse pass, the~$I^{(R)}$ and~$I^{(C)}$ indices switch roles.
This cycle constitutes one full iteration of the TTOpt method.
Note that the procedure described above assumes fixed ranks; however, as mentioned earlier, maxvol can be used to select more rows to form rectangular matrices.
Consequently, an adaptive version of the algorithm that allows for increasing ranks can be considered.

These iterations are repeated until either the budget for black-box calls is exhausted or the algorithm converges, \ie, no index changes during a cycle.
Ultimately, the quasi-optimal black-box element is returned along with its corresponding indices as the optimum among all elements queried during the iterations.
Note that TTOpt exploits the low-rank TT-structure of $\tens{E}$, requiring only $\mathcal{O}(dR^2 \max N_i)$ operations per iteration, where $R$ is the average TT-rank.

\subsubsection{Method PROTES}

The PROTES (PRobabilistic Optimizer with TEnsor Sampling) approach reformulates optimization as a probability density estimation problem.
It constructs a parametric distribution $p_{\boldsymbol{\theta}}(\vect{n})$ over the search space $\Phi$ using a TT-tensor, wherein high values (probabilities) correspond to low energies.
This probability distribution is modeled as a low-parameter black box satisfying two conditions:
we must be able to sample from it according to the specified probability (\ie, the probability of a sample $\vx_0\in \Phi$ occurring is proportional to the black box value $\distp_\theta(\vx)$ at that point, where~$\theta$ represents the set of all black-box parameters);
and we must be able to evaluate the black box (\ie, the unnormalized probability) at an arbitrary point.
Due to its low-rank property and the availability of necessary operations, the TT-decomposition is well-suited for this role.

The general structure of the PROTES algorithm, adapted here for the optimization of atomic clusters, is presented in Figure~\ref{fig:scheme_protes}.
The method is designed to find a quasi-optimal point of a given discrete multivariate function \(E\).
Starting with an initial tensor \(\tens{E}_{0} \) in the TT-format, the algorithm iteratively refines the solution using stochastic gradient descent (SGD) with a specified learning rate \(\eta\).
The parameters \(K\) and \(k\) control the number of sampled points and the selection of the top-performing candidates, respectively.
The main steps of the algorithm are described below.

\textbf{Sample.}
At each iteration, until a stopping criterion is met (such as exhausting the computational budget), we sample~\cite{dolgov2020approximation} a set of \(K\) discrete candidates $\{\vect{n}^{(i)}\}_{i=1}^K$ from the probability distribution $p_{\boldsymbol{\theta}}(\vect{n})$ represented as a TT-decomposition.

In our experiments, to ensure non-negativity, we model the probability~$p_{\boldsymbol{\theta}}$ as being proportional to the square of the tensor.
This implies that the probability of sampling the configuration $\{n_1,n_2,\ldots,n_d\}$ is:
\begin{equation}\label{eq:dens_in_TT}
    p_{\boldsymbol{\theta}}(\{n_1,n_2,\ldots,n_d\})=\frac1Z \mathcal P^2[n_1,n_2,\ldots,n_d],
\end{equation}
where $\mathcal P$ is the tensor whose cores we optimize, and $Z$ is the normalizing constant:
$$
Z=\sum_{\xi_1,\xi_2,\ldots,\xi_d}\mathcal P^2[\xi_1,\xi_2,\ldots,\xi_d].
$$
Note that this sampling procedure, described in detail in~\cite{dolgov2020approximation}, is computationally efficient and relies on sequential sampling of each coordinate based on marginal distributions.
Specifically, for any multivariate discrete probability~$p$, we first calculate the marginal (independent of other components) probability distribution~$p_1$ for the first component:
$$
p_1(n_1)=\sum_{n_2,\ldots, n_d}{p(n_1,n_2,\ldots, n_d)}.
$$
We then sample from this distribution, which is straightforward since $p_1(n_1)$ is simply a vector of probabilities.
Let $\hat n_1$ be the value obtained.
In the next step, we calculate the conditional marginal distribution~$p_2$ of the second component, given that the first index has the value sampled in the previous step:
\begin{multline}
p_2(n_2)=
\frac{\sum_{n_3,\ldots, n_d}{p(\hat n_1,n_2,n_3\ldots, n_d)}}{p_1(\hat n_1)}=\\
\frac{\sum_{n_3,\ldots, n_d}{p(\hat n_1,n_2,n_3\ldots, n_d)}}%
  {\sum_{n_2,n_3\ldots, n_d}{p(\hat n_1,n_2,n_3\ldots, n_d)}}
  .
\end{multline}
Proceeding recursively, at the $k^{\text{th}}$ step, we have the sampled components $\{\hat n_1, \hat n_2, \ldots, \hat n_{k-1}\}$ and calculate the conditional marginal distribution~$p_k$ of the $k^{\text{th}}$ component:
$$
p_k(n_k)=
\frac%
{\sum_{n_{k+1},\ldots, n_d}{p(\hat n_1, \ldots, \hat n_{k-1}, n_k, n_{k+1}\ldots, n_d)}}%
{\sum_{n_k,\ldots, n_d}{p(\hat n_1, \ldots, \hat n_{k-1}, n_k, n_{k+1}\ldots, n_d)}}.
$$
The primary challenge lies in accurately computing these multidimensional sums.
The TT-decomposition addresses this effectively by allowing the efficient calculation of such sums.
Consider a probability density in the form~\eqref{eq:dens_in_TT}, where the TT decomposition is represented as in~\eqref{eq:tt_repr}.
First, we require the TT-decomposition of the tensor square.
This is achieved by taking the Kronecker product of each core slice.
Specifically, each core~$\mathcal G_i$ is replaced by a core with a higher rank~$\mathcal G'_i$ as follows:
$$
\mathcal G'_i(\overline{jk},\,n,\,\overline{lm})=
\mathcal G_i(j,\,n,\,l)\mathcal G_i(k,\,n,\,m),
$$
where $\overline{jk}$ and $\overline{lm}$ are super-indices, defined in~\eqref{eq:def_super-indices}.

Subsequently, the required sums can be expressed as:
\begin{multline*}
    \sum_{n_{k+1},\ldots, n_d}\mathcal P^2[\hat n_1, \hat n_2, \ldots, \hat n_{k-1}, n_k, n_{k+1}\ldots, n_d]=\\=
\sum_{r_0=1}^{R_0}
\sum_{r_1=1}^{R_1}
\cdots
\sum_{r_{d}=1}^{R_{d}}
    \tens{G}'_1 [r_0,\hat n_1 , r_1]
    \tens{G}'_2 [r_1, \hat n_1, r_2]
    \times
    \cdots
    \times
    \\
    \times
    \tens{G}'_{k-1} [r_{k-2}, \hat n_{k-1}, r_{k-1}]
    \tens{G}'_{k} [r_{k-1}, n_{k}, r_{k}]
    \times
    \\
    \times
    \left( 
    \sum_{n_{k+1}=1}^{N_{k+1}}
    \tens{G}'_{k+1} [r_{k}, \hat n_{k+1}, r_{k+1}]
    \right)
    \times
    \cdots
    \times
    \\
    \times
    \left( 
    \sum_{n_{d}=1}^{N_{d}}
    \tens{G}'_{d} [r_{d-1}, \hat n_{d}, r_{d}]
    \right).
\end{multline*}
Note that the computational complexity of evaluating such convolutions (considering only multiplications) is equivalent to that of calculating a single tensor element.

\textbf{Approve.}
We evaluate the $K$ multi-indices proposed during sampling to check for compliance with physical constraints and select the $\tilde{K}$ most relevant candidates.

\textbf{Compute.}
Energies $\{e_i = \tens{E}[\vect{n}^{(i)}]\}_{i=1}^{\tilde{K}}$ are computed for each of the proposed candidates.
This evaluation is typically executed in parallel to efficiently utilize computational resources.

\textbf{Select.}
We select the \(k\) best points based on their energy values.
For a minimization problem, these correspond to the points with the lowest values of \(e_i\): 
\[
\mathcal{J} = \{ j : e_j \leq e_{(k)} \}, \quad |\mathcal{J}| = k
\]

\textbf{Update.}
We adjust the parameters $\boldsymbol{\theta}$ (\ie, the core elements of the tensor $\mathcal P$) to increase the likelihood of the indices in $\mathcal{J}$.
Specifically, we maximize the likelihood by minimizing the loss~$\mathcal{L}(\boldsymbol{\theta})$:
$$
\mathcal{L}(\boldsymbol{\theta}) = -\sum_{j \in \mathcal{J}} \log p_{\boldsymbol{\theta}}(\vect{n}^{(j)})
$$
based on which we perform a gradient descent step:
$$
\boldsymbol{\theta} \leftarrow \boldsymbol{\theta} + \eta \nabla_{\boldsymbol{\theta}} \mathcal{L}(\boldsymbol{\theta}).
$$
Thus, PROTES effectively balances exploration and exploitation by sampling a diverse set of points while progressively focusing on regions of the search space that yield better function values.

Note that fundamentally, TTOpt exploits algebraic structures within $\tens{E}$, while PROTES leverages the information geometry \cite{Amari2000Methods, e22101100} of the solution space.
This makes the methods complementary: TTOpt is preferable when $\tens{E}$ admits an accurate low-rank approximation, whereas PROTES excels on rugged energy surfaces with sparse minima.
Moreover, both methods overcome the curse of dimensionality inherent in Eq.~\eqref{eq:opt_def} by leveraging the low-rank TT-decomposition.

\subsubsection{Physically-constrained initialization with PROTES}

As is evident from the description of the encodings above (see Seq.~\ref{seq:grids}), the algorithm allows for efficient operation in both directions: decoding the floating-point positions of particles from the integer indices output by the tensor, and encoding, with a certain degree of accuracy, the floating-point parameter values into integer indices compatible with the tensor mechanisms.

We utilize this property to initialize the tensor.
Before describing the initialization mechanism in detail, it is important to note that converting floating-point coordinates to integers may utilize a lossy representation.
That is, if a floating-point value~$x$ is encoded to a corresponding integer and then decoded via $x'=i2f\bigl( f2i(x) \bigr)$, the values~$x$ and~$x'$ coincide only approximately; the proximity is determined by the bit depth and the physical restrictions imposed by the specific encoding.
However, this approximation is acceptable for initialization purposes, as starting with a reasonably feasible configuration is more critical than exact reproduction.

Assume we select~$r$ to be approximately equal or slightly larger than the target rank and consider~$r$ configurations (which can be obtained using PyXtal or similar libraries).
We encode each configuration using the selected encoding~$X_0=\{x_0,\,\ldots,\,x_{d-1}\}\in\mathbb N^d$.
Next, we construct a rank-1 TT-tensor that equals one when its indices coincide with the components of this vector~$X_0$, and zero otherwise.
The cores~$\{\mathcal G_i\}_{i=1}^d$ of such a decomposition are constructed explicitly as follows:
$$
G_i[1, n, 1]=\begin{cases}
    1, &n=x_i,\\
    0, &\text{otherwise}.
\end{cases}
$$

We obtain~$r$ such uniform-rank TT-tensors.
These tensors are then summed (via an efficient operation in the TT-format), the result is rounded if the specified rank is less than~$r$, and a small amount of noise is added to the cores.
This noise addition prevents zero-valued entries and ensures the resulting tensor can generate a sufficiently diverse set of configurations beyond the initial~$r$ examples.


\subsection{Potentials under consideration}
    \subsubsection{The Lennard-Jones potential}

The Lennard-Jones (LJ) potential, \( E_{LJ} \), represents one of the simplest and most fundamental models for describing the interaction between neutral atoms or molecules in molecular dynamics simulations and various other applications.
This potential is defined by the following expression:
\begin{align} \label{eq:LJ}
E_{LJ}(r) = 4\epsilon \left[ \left( \frac{\sigma}{r} \right)^{12} - \left( \frac{\sigma}{r} \right)^{6} \right]
\end{align}
The parameters governing the LJ potential are defined as follows: \(\epsilon\) denotes the depth of the potential well, representing the strength of attraction between two particles at the equilibrium distance; specifically, a larger value of \(\epsilon\) corresponds to stronger attraction.
The parameter \(\sigma\) indicates the finite distance at which the inter-particle potential equals zero, effectively marking the transition point between attractive and repulsive interactions as particles approach one another.
Finally, \(r\) denotes the distance between the centers of the two particles.

\subsubsection{Moment Tensor Potential and its fitting}

In this study, we also employ Moment Tensor Potential (MTP)~\cite{shapeev2016_mtp}.
The MTP framework serves as a local interaction model wherein the total energy is calculated as the sum of contributions from each $i$-th atom within the C20 structure:
\begin{equation} \label{eq:MTP}
    E_{\rm MTP} = \sum \limits_i \sum \limits_{\alpha} \xi_{\alpha} B_{\alpha} ({\bf \mathfrak{n}}_i).
\end{equation}
In this equation, $\xi_{\alpha}$ represents the linear MTP parameters, while $B_{\alpha}({\bf \mathfrak{n}}_i)$ denotes the basis functions dependent on the local neighborhood, ${\bf \mathfrak{n}}_i$, of the $i$-th atom.
Each atomic neighborhood ${\bf \mathfrak{n}}_i$ comprises the central $i$-th atom at position ${\bf r}_i$ and all neighboring $j$-th atoms at positions ${\bf r}_j$ that satisfy the condition $r = |{\bf r}_{ij}| = |{\bf r}_j - {\bf r}_i| \leq R_{\rm cut}$, where $R_{\rm cut}$ refers to the cutoff radius.
Consequently, only interactions between the central atom and neighbors situated within a sphere of radius $R_{\rm cut}$ are considered.
The basis functions $B_{\alpha}$ are constructed via the contraction of the Moment Tensor Descriptors:
\begin{align} \label{Descriptor}
    M_{\mu, \nu}({\bf \mathfrak{n}}_i) = \sum \limits_j f_{\mu}(|{\bf r}_{ij}|) {\bf r}_{ij}^{\otimes \nu},
\end{align}
where the radial component, $f_{\mu}(|{\bf r}_{ij}|)$, takes the following form:
\begin{equation} \label{RadialPart}
    f_{\mu}(|{\bf r}_{ij}|) = \sum \limits_{\beta} c_{\mu, \beta} T^{\beta}(|{\bf r}_{ij}|) (R_{\rm cut} - |{\bf r}_{ij}|)^2.
\end{equation}
Here, $\mu$ is the index of the radial function, $T^{\beta}$ represents the Chebyshev polynomial of order $\beta$, and $c_{\mu, \beta}$ are the radial MTP parameters.
The symbol ``$\otimes$'' denotes the outer product of the vectors ${\bf r}_{ij}$, applied $\nu$ times; thus, the angular component ${\bf r}_{ij}^{\otimes \nu}$ constitutes a tensor of rank $\nu$.
The indices $\mu$ and $\nu$ determine the level of the Moment Tensor Descriptor, given by ${\rm lev} M_{\mu, \nu} = 2 + 4 \mu + \nu$.
Correspondingly, the level of the MTP basis function is defined as ${\rm lev} B_{\alpha} = \sum \limits_m (2 + 4 \mu_m + \nu_m)$, which represents the sum of the levels of the Moment Tensor Descriptors that form a scalar contraction.
To establish a specific functional form for the MTP, we select a maximum potential level, ${\rm lev}_{\rm MTP}$, and encompass only those basis functions satisfying the condition ${\rm lev} B_{\alpha} \leq {\rm lev}_{\rm MTP}$.

The potential parameters ${\bm \theta} = \{\xi_{\alpha}, c_{\mu, \beta}\}$ are determined by minimizing the following objective function:
\begin{align*} 
\displaystyle
\min\limits_{\bm \theta} \quad
\sum \limits_{k=1}^K \left[
	\left(E_k({\bm \theta}) - E^{\rm DFT}_k \right)^2
	+
	w_{\rm f} \sum_{i} \left| {\bm f}_{k,i}({\bm \theta}) - {\bm f}^{\rm DFT}_{k,i} \right|^2 \right]
\end{align*}
where $K$ denotes the total number of structures within the training set.
The terms $E^{\rm DFT}_k$ and ${\bm f}^{\rm DFT}_{k,i}$ represent the energies and forces, respectively, calculated for the $k$-th structure using density functional theory (DFT).
Similarly, $E_k({\bm \theta})$ and ${\bm f}_{k,i}({\bm \theta})$ denote the energies and forces predicted by the MTP as functions of the parameters ${\bf \theta}$.
Finally, $w_{\rm f}$ is a non-negative weighting factor that quantifies the relative importance of fitting the MTP forces to the DFT forces within the objective function.
To automatically generate the training set for MTP fitting, we employed an active learning (AL) algorithm briefly described below.

Assume that we find the optimal parameters ${\bar{\bm \theta}}$ after fitting the MTP on some initial training set. Then, we compose the following matrix
\[
\mathsf{B}=\left(\begin{matrix}
\frac{\partial E_1\left( {\bar{\bm \theta}} \right)}{\partial \theta_1} & \ldots & \frac{\partial E_1\left( {\bar{\bm \theta}} \right)}{\partial \theta_m} \\
\vdots & \ddots & \vdots \\
\frac{\partial E_K\left( {\bar{\bm \theta}} \right)}{\partial \theta_1} & \ldots & \frac{\partial E_K\left( {\bar{\bm \theta}} \right)}{\partial \theta_m} \\
\end{matrix}\right),
\]
where each row corresponds to a particular structure. Next, we construct a subset of structures yielding the most linearly independent rows (physically it means geometrically different structures) in $\mathsf{B}$. This is equivalent to finding a square $m \times m$ submatrix $\mathsf{A}$ of the matrix $\mathsf{B}$ of maximum volume. To that end, we use the so-called maxvol algorithm \cite{zamarashkin2010-maxvol}. To determine whether a given structure $\bm x^*$ obtained during structural relaxation (geometry optimization) is representative, we calculate the extrapolation grade $\gamma(\bm x^*)$ defined as
\begin{equation} \label{Grade}
\begin{array}{c}
\displaystyle
\gamma(\bm x^*) = \max_{1 \leq j \leq m} (|c_j|), ~\rm{where}
\\
\displaystyle
{\bm c} = \left( \dfrac{\partial E}{\partial \theta_1} (\bar{\theta}, \bm x^*) \ldots \dfrac{\partial E}{\partial \theta_m} (\bar{\theta}, \bm x^*) \right) \mathsf{A}^{-1}.
\end{array}
\end{equation}
This grade defines the maximal factor by which the determinant $|{\rm det(\mathsf{A})}|$ can be increased if ${\bm x^*}$ is added to the training set. Thus, if the structure $\bm x^*$ is a candidate for adding to the training set then $\gamma_{\rm th} \leq \gamma(\bm x^*) \leq \Gamma_{\rm th}$, where $\gamma_{\rm th} = 2.1$ is an adjustable threshold parameter which controls the value of permissible extrapolation. Once the extrapolation grade exceeds $\Gamma_{\rm th} = 10$ the relaxation is stopped. Such values of the thresholds $\gamma_{\rm th}$ and $\Gamma_{\rm th}$ were chosen based on the previous benchmarks~\cite{Podryabinkin2017,novikov2020_mlip_2}. Once the relaxation is stopped, we calculate the matrix with the derivatives of energies for the candidate structures with respect to the MTP parameters and combine it with the matrix $\mathsf{A}$. Using the maxvol algorithm, we select the rows corresponding to different structures in the combined matrix, calculate the new matrix $\mathsf{A}$, and incorporate all the novel different structures into the training set. Next, we conduct DFT calculations and obtain DFT energies and forces for the novel structures. Finally, we re-fit the potential. This procedure repeats until all structural optimizations can run without preselecting structures, i.e., without the occurrence of the candidate structures. 
\section{RESULTS AND DISCUSSION}
    \label{sec:results}
    \subsection{Experiments with Lennard-Jones clusters}
\label{sec:results_lj}

We validated our TT-based optimization framework using the well-studied LJ potential~\eqref{eq:LJ}, adopting standard parameters $\epsilon = 1$ and $\sigma = 1$ in reduced units.
Our analysis reveals a fundamental trade-off between computational efficiency and reliability that is critically dependent on the size of the cluster, the encoding scheme and the initialization strategy.
We selected clusters ranging from 5 to 45 atoms to encompass different optimization regimes~\cite{Freitas_Gustavo2022-va, GAS, WalesLJ38}.
Small clusters (5--15 atoms) establish baseline benchmarks, medium clusters (16--30 atoms) exhibit increasingly rugged energy landscapes, and larger systems (31--45 atoms) are used to assess scalability.
Note that the LJ$_{38}$ cluster constitutes a particularly challenging benchmark due to its double-funnel energy landscape, which frequently entraps optimization algorithms.

Performance was evaluated using two complementary metrics: PROTES Calls (PC), defined as the number of evaluations of the energy function, representing search efficiency; and the success rate ($SR_t$), defined as the percentage of independent runs that locate the global minimum, representing reliability.
Comprehensive results, along with computational experiment parameters and their sensitivity analysis, are provided in the Supplementary Material \ref{sec:supp}.

\subsubsection{The efficiency-reliability trade-off}

Table~\ref{table:consolidated_performance} illustrates a distinguishable pattern in which no single configuration yields optimal performance across the entire range of cluster sizes.
Consequently, a dichotomy emerges, necessitating a choice between high-efficiency and high-reliability operational modes.

\begin{table}[ht]
    \centering
    \small
    \caption{Comparative performance of the PROTES optimizer for selected Lennard-Jones clusters. The table shows the number of PROTES calls  (PC, global energy evaluations) and success rate ($SR_t$) for different combinations of encoding schemes and initialization strategies. Encoding schemes: SR = Simple Relative, CR = Constrained Relative. Initialization: Agn =  Agnostic, PhC = Physically-Constrained. For results using the TTOpt algorithm with direct encoding, see corresponding Table in the Appendix.}
    \label{table:consolidated_performance}
    \begin{tabular}{|c|l|c|r|r|r|}
    \hline
    $N_a$ & \multicolumn{1}{c|}{Method} & Init. & \multicolumn{1}{c|}{PC} & \multicolumn{1}{c|}{$SR_t$ (\%)} & \multicolumn{1}{c|}{Energy (eV)} \\
    \hline
    \multirow{4}{*}{13} & PROTES + SR & Agn & \textbf{104} & 100 & -44.33 \\
                        & PROTES + SR & PhC & 3,200 & 100 & -44.33 \\
                        & PROTES + CR & Agn & 25,400 & 100 & -44.33 \\
                        & PROTES + CR & PhC & 1,690 & 100 & -44.33 \\
    \hline
    \multirow{4}{*}{18} & PROTES + SR & Agn & 6,655 & 100 & -66.53 \\
                        & PROTES + SR & PhC & \textbf{157,000} & 100 & -66.53 \\
                        & PROTES + CR & Agn & 7,125 & 100 & -66.53 \\
                        & PROTES + CR & PhC & 269,000 & 100 & -66.53 \\
    \hline
    \multirow{4}{*}{26} & PROTES + SR & Agn & 2,147 & 100 & -108.32 \\
                        & PROTES + SR & PhC & \textbf{234,000} & 100 & -108.32 \\
                        & PROTES + CR & Agn & 17,452 & 100 & -108.32 \\
                        & PROTES + CR & PhC & 1,040,000 & 100 & -108.32 \\
    \hline
    \multirow{4}{*}{33} & PROTES + SR & Agn & 22,100 & 41 & -161.92 \\
                        & PROTES + SR & PhC & 22,100 & 41 & -161.92 \\
                        & PROTES + CR & Agn & 23,400 & 25 & -161.92 \\
                        & PROTES + CR & PhC & \textbf{22,500} & 25 & -161.92 \\
    \hline
    \multirow{2}{*}{38} & PROTES + SR & PhC & \textbf{48,900} & 16 & -173.93 \\
                        & PROTES + CR & PhC & 54,500 & 33 & -173.93 \\
    \hline
    \multirow{2}{*}{45} & PROTES + SR & PhC & 55,600 & 0 & -- \\
                        & PROTES + CR & PhC & \textbf{52,300} & 33 & -244.33 \\
    \hline
    \end{tabular}
\end{table}

For clusters containing up to 26 atoms, agnostic initialization paired with Simple Relative encoding demonstrates superior efficiency, achieving 100\% success rates with significantly reduced computational budgets (e.g., $2,147$ evaluations for LJ$_{26}$).
This corresponds to an ``efficiency-optimized'' configuration, wherein the TT-model rapidly identifies promising regions in the absence of physical bias.

However, for larger clusters ($N_a \geq 30$), the robustness of this encoding configuration diminishes significantly.
In this regime, physically-constrained initialization becomes essential, effectively transitioning the method into a ``reliability-optimized'' mode.
Although this approach requires more evaluations, it maintains non-zero success rates in cases where agnostic initialization fails to converge.
Notably, for LJ$_{45}$, only physically-constrained initialization used in conjunction with Constrained Relative encoding achieves a successful outcome (33\% SR, 52,300 evaluations).

\subsubsection{Encoding strategy analysis}

The selection of an encoding scheme further modulates this trade-off.
Specifically, Simple Relative encoding generally offers enhanced efficiency for smaller clusters and practically facilitates accelerated convergence when successful.
In contrast, Constrained Relative encoding provides superior reliability for challenging systems, as evidenced by success rates on LJ$_{38}$ (33\% vs 16\%) and LJ$_{45}$ (33\% vs 0\%) that exceed those of the simple variant.

The TTOpt algorithm utilizing direct encoding (see Table~\ref{table:baselines} and Table~\ref{tab:ttopt_results} in supplementary material) occupies a distinct position within this design space, characterized by maximum reliability (100\% $SR_t$ across tested clusters) but significantly higher computational costs.
This demonstrates that while advanced encodings dramatically improve efficiency, the underlying TT-approach remains fundamentally robust.

\subsubsection{Performance on challenging landscapes}

\begin{figure}[t!]
    \centering
    \includegraphics[width=0.95\linewidth]{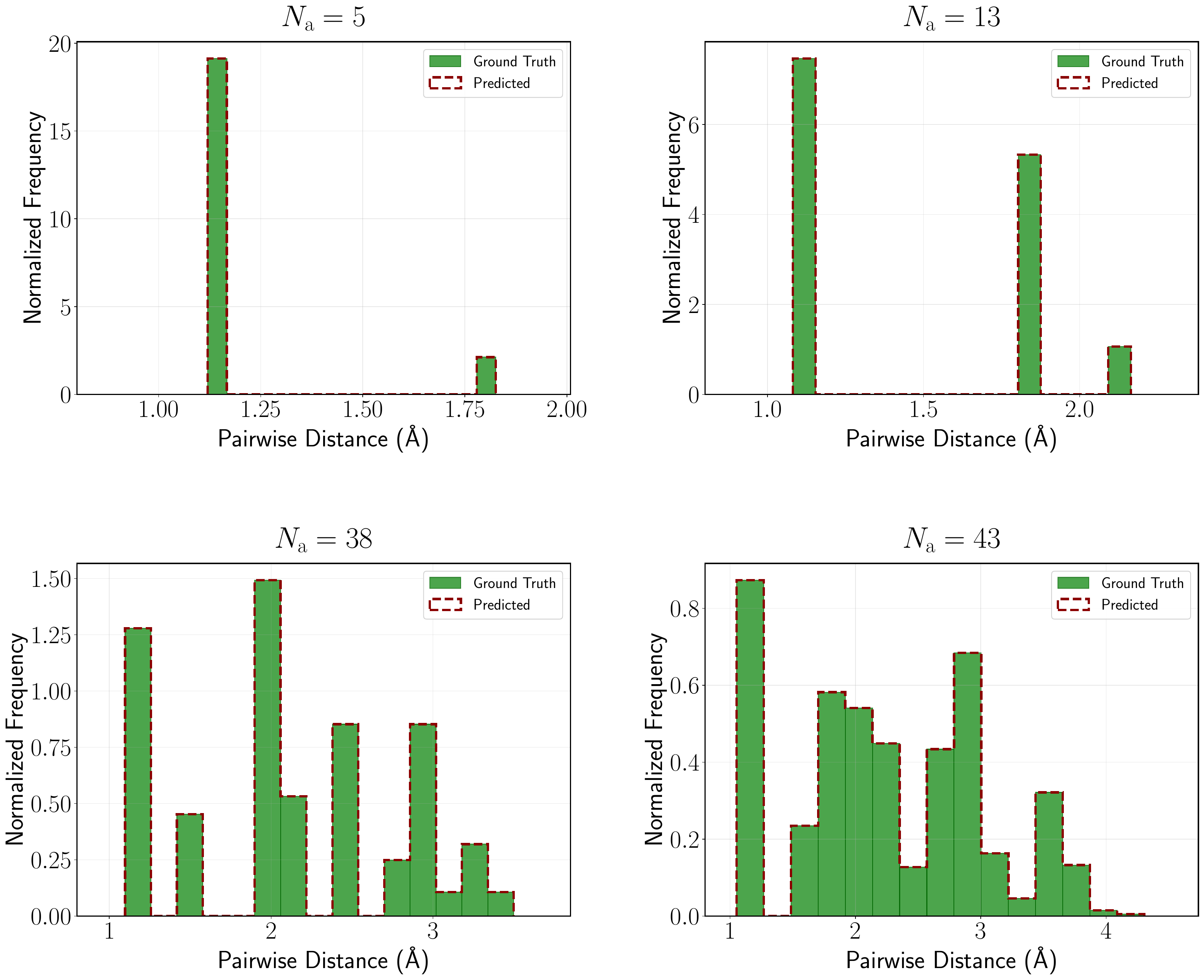}
    \caption{
        Validation of optimized LJ cluster structures through pairwise distance analysis.
        The ground truth distributions (green) are closely matched by our method's predictions (red dashed) for clusters containing 5, 13, 38, and 43 atoms.
        The accurate reproduction of distance histograms confirms the structural correctness of the identified minima.
    }
    \label{fig:lj_pairwise_distributions}
\end{figure}

The results obtained for LJ$_{38}$ warrant particular attention.
Our method achieves modest yet meaningful success rates (16-33\%) on this notoriously difficult system, confirming the framework's capability to navigate complex, multi-funnel energy landscapes.
The accurate reproduction of the LJ$_{38}$ global minimum structure is substantiated by its pairwise distance distribution (Figure~\ref{fig:lj_pairwise_distributions}).
Although scope for further optimization exists, these results represent a significant achievement considering the general-purpose nature of the method and its minimal reliance on problem-specific tuning.

\subsubsection{Comparative analysis with established methods}

Table~\ref{table:baselines} contextualizes our findings within the broader optimization literature.
A distinguishing feature of our approach is the ``single-hop'' refinement strategy.
Whereas methods such as USPEX~\cite{USPEX2009} necessitate numerous local relaxations, our method typically requires only a single relaxation, potentially offering substantial advantages for computationally expensive quantum-mechanical potentials.
The probabilistic PROTES method coupled with advanced encodings demonstrates superior efficiency for tractable problems, while the algebraic TTOpt provides guaranteed reliability at a higher computational cost.
This duality represents a key strength of the proposed framework, enabling users to select the appropriate balance for their specific application.

\begin{table}[ht]
\small
\caption{Comparison of optimization efficiency for Lennard-Jones clusters. Our TT-based methods achieve competitive performance in terms of global calls (C) with ``single hop'' local refinement (L = 1), while established methods require multiple local optimizations. Best results for each cluster size are highlighted. Encoding names are abbreviated as follows: SR = Simple Relative, CR = Constrained Relative.}
\label{table:baselines}
\centering
\begin{tabular}{|c|c|c|c|c|}
\hline
$N_{a}$ & M & C & L & Energy \\
\hline
13 & PSO~\cite{MAI2020100710} & - & - & -44.2 \\
13 & TTopt & 10,729 & 1 & -44.3 \\
13 & PROTES + SR & \textbf{104} & \textbf{1} & -44.3 \\
13 & PROTES + CR & 25,385 & 1 & -44.3 \\
\hline
18 & PSO~\cite{MAI2020100710} & - & - & -66.2 \\
18 & TTopt & 29,451 & 1 & -66.53 \\
18 & PROTES + SR & \textbf{6,655} & \textbf{1} & -66.53 \\
18 & PROTES + CR & 7,125 & 1 & -66.53 \\
\hline
26 & USPEX~\cite{USPEX2009} & 34,200 & 56 & -108.3156 \\
26 & Locatelli2003~\cite{locatelli2003} & 36,514 & 267 & -108.3156 \\
26 & Minima Hopping~\cite{Goedecker2004} & 50,610 & 96 & -108.3156 \\
26 & GAS~\cite{GAS} & 12,655 & 34 & -108.3156 \\
26 & TTopt & 19,544 & 1 & -108.3156 \\
26 & PROTES + SR & \textbf{2,147} & \textbf{1} & -108.3156 \\
26 & PROTES + CR & 17,452 & 1 & -108.3156 \\
\hline
38 & GAS~\cite{GAS} & 278,828 & 752 & -173.9284 \\
38 & PROTES + SR & \textbf{42,415} & 1 & -173.9284 \\
\hline
\end{tabular}
\end{table}

\subsection{Experiment with C20 and MTP}

As a second, more realistic case study, we consider carbon clusters consisting of 20 atoms, and to model this system, we employed a level-24 MTP, which includes 913 parameters ${\bf \theta}$.
We used a cutoff radius of $R_{\rm cut} = 5 ~\angstrom$.
To generate an initial training set of $\text{C}_{20}$ configurations, we utilized the random-symmetric atomic cluster generator implemented in the PyXtal package~\cite{fredericks2021pyxtal}, yielding a set of 200 clusters.
We subsequently performed DFT calculations on these clusters, retaining only those for which electronic convergence was achieved.
The resulting initial training set consisted of 174 converged clusters, upon which an initial MTP was fitted.

Next, we relaxed the geometries of these clusters using the initial MTP and iteratively expanded the training set via an active learning (AL) strategy based on the MaxVol algorithm (see, e.g.,~\cite{novikov2020_mlip_2}).
This procedure yielded the final fitted MTP and a training set comprising 2737 configurations.
All DFT calculations were performed using the VASP package.
We employed the Perdew-Burke-Ernzerhof generalized gradient approximation (PBE-GGA) as the exchange–correlation functional, with a plane-wave energy cutoff of 550 eV.
Each cluster was placed in a simulation box with sufficient vacuum to prevent periodic interactions; consequently, the Brillouin zone was sampled at the $\Gamma$-point.

\begin{figure}[h!]
    \centering
    \includegraphics[width=1.0\linewidth]{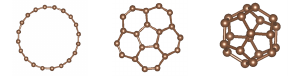}
    \caption{
        Structures of the identified $\text{C}_{20}$ atomic clusters: monocyclic cap, fullerene cage structure, and buckyball.
    }
    \label{fig:C20}
\end{figure}

Upon relaxing the $\text{C}_{20}$ clusters using the final MTP, we identified three distinct local minima.
The first corresponds to a monocyclic cap with an energy of $-8.141$ eV/atom, the second is a fullerene cage structure ($-8.041$ eV/atom), and the third is a buckyball with an energy of $-8.039$ eV/atom (see Fig.~\ref{fig:C20}).
The monocyclic cap and the fullerene cage have been previously reported, for example, in~\cite{deaven1995_c_clusters} using a genetic algorithm, while the buckyball structure was identified in~\cite{cai2004_brenner} using the Brenner potential combined with a global optimization algorithm. We note that the buckyball does not correspond to the lowest energy, as the MTP fitting error was approximately 70 meV/atom. Nevertheless, our primary objective was to test the proposed global optimization algorithms rather than to accurately fit an MTP for the correct reproduction of the lowest energy.

To validate the performance of the proposed global optimization algorithms, we applied the trained MTP to this system.

Both TTopt and PROTES algorithms identified the monocyclic cap structure as the global minimum for C$_{20}$ clusters. Across multiple independent runs, the mean optimized total cluster energy was $-162.829 \pm 0.0010 \text{ eV}$. This corresponds to a potential energy per atom of $-8141.45\text{ meV/atom}$, with a standard deviation of $0.05 \text{ meV/atom}$ across runs. This result aligns precisely with the global minimum identified during the MTP fitting procedure.

\section{CONCLUSION}
    \label{sec:conclusions}
    In this work, we introduce a novel framework that integrates tensor train decomposition with physically-constrained encoding schemes for the global optimization of atomic clusters.
Although the application to Lennard-Jones clusters containing up to 45 atoms showcases the method's capabilities, the primary contribution lies in the methodological framework itself.
This framework establishes a new paradigm for navigating high-dimensional potential energy surfaces and can be integrated with established optimization techniques~\cite{USPEX2009, Wales1997, Wales2003, CSPBench}.

A key innovation of our approach is the direct encoding of physical constraints (such as plausible interatomic distances and angular relationships) into the tensor decomposition process.
Crucially, this physically constrained formulation dramatically reduces the effective dimensionality of the search space, thereby preserving mathematical rigor while bridging the gap between abstract optimization strategies and chemical intuition.
Our comprehensive analysis reveals a fundamental trade-off between efficiency and reliability that governs the framework's performance.
For small- to medium-sized clusters ($N_a \leq 26$), the probabilistic PROTES method, utilizing Simple Relative encoding and agnostic initialization, achieves remarkable efficiency, often locating global minima with orders of magnitude fewer evaluations than established benchmarks.
Conversely, for larger and more challenging systems ($N_a \geq 30$), physically-constrained initialization proves essential for robustness, while the algebraic TTOpt method offers guaranteed reliability at a higher computational cost.

This duality constitutes a key strength of the framework, enabling users to select an appropriate balance for their specific applications.
The successful optimization of LJ$_{38}$, which is a benchmark system known for its challenging double-funnel landscape~\cite{WalesLJ38}, demonstrates the method's capacity to navigate complex energy surfaces characteristic of real materials.
Furthermore, the application to carbon clusters using a machine-learning Moment Tensor Potential confirms the framework's practical utility beyond LJ potentials, as it achieved geometries consistent with quantum-accurate simulations.
The framework's ``single-hop'' refinement strategy, which typically requires only one local relaxation, renders it particularly valuable for quantum-mechanical potentials where force calculations dominate the computational cost.

Future work will focus on developing adaptive schemes that automatically tune discretization parameters and initialization strategies during optimization, as well as hybrid approaches that combine our tensor methods with USPEX~\cite{USPEX2009,USPEX2011}, Basin Hopping~\cite{Wales1997,Wales2003}, and other established frameworks~\cite{CSPBench}.
The capacity to optimize structures using the same algorithm in diverse potential classes represents a significant step toward automated materials discovery pipelines.
By providing a mathematically rigorous yet physically intuitive approach to high-dimensional optimization, this work opens new pathways for analyzing complex energy landscapes in computational materials science.
    
\section*{Acknowledgments}
The work was supported by the grant for research centers in the field of AI provided by the Ministry of Economic Development of the Russian Federation in accordance with the agreement 000000C313925P4F0002 and the agreement with Skoltech №139-10-2025-033 from 20.06.2025.

\bibliography{references.bib}

@inproceedings{TTopt_NEURIPS2022,
 author = {Sozykin, Konstantin and Chertkov, Andrei and Schutski, Roman and Phan, Anh-Huy and Cichocki, Andrzej S and Oseledets, Ivan},
 booktitle = {Advances in Neural Information Processing Systems},
 editor = {S. Koyejo and S. Mohamed and A. Agarwal and D. Belgrave and K. Cho and A. Oh},
 pages = {26052--26065},
 publisher = {Curran Associates, Inc.},
 title = {{TTOpt}: A Maximum Volume Quantized Tensor Train-based Optimization and its Application to Reinforcement Learning},
 url = {https://proceedings.neurips.cc/paper_files/paper/2022/file/a730abbcd6cf4a371ca9545db5922442-Paper-Conference.pdf},
 volume = {35},
 year = {2022}
}

@article{fredericks2021pyxtal,
  title={PyXtal: A Python library for crystal structure generation and symmetry analysis},
  author={Fredericks, Scott and Parrish, Kevin and Sayre, Dean and Zhu, Qiang},
  journal={Computer Physics Communications},
  volume={261},
  pages={107810},
  year={2021},
  publisher={Elsevier}
}

@article{MAI2020100710,
title = {Optimization of Lennard-Jones clusters by particle swarm optimization with quasi-physical strategy},
journal = {Swarm and Evolutionary Computation},
volume = {57},
pages = {100710},
year = {2020},
issn = {2210-6502},
doi = {https://doi.org/10.1016/j.swevo.2020.100710},
url = {https://www.sciencedirect.com/science/article/pii/S2210650220303631},
author = {Guizhen Mai and Yinghan Hong and Shen Fu and Yingqing Lin and Zhifeng Hao and Han Huang and Yuanhao Zhu},
}

@article{oseledets2011tensor,
  title={Tensor-train decomposition},
  author={Oseledets, Ivan V},
  journal={SIAM Journal on Scientific Computing},
  volume={33},
  number={5},
  pages={2295--2317},
  year={2011},
  publisher={SIAM}
}

@inproceedings{NEURIPS2023_02895786,
 author = {Batsheva, Anastasiia and Chertkov, Andrei and Ryzhakov, Gleb and Oseledets, Ivan},
 booktitle = {Advances in Neural Information Processing Systems},
 editor = {A. Oh and T. Naumann and A. Globerson and K. Saenko and M. Hardt and S. Levine},
 pages = {808--823},
 publisher = {Curran Associates, Inc.},
 title = {PROTES: Probabilistic Optimization with Tensor Sampling},
 url = {https://proceedings.neurips.cc/paper_files/paper/2023/file/028957869e560af14243ac37663a471e-Paper-Conference.pdf},
 volume = {36},
 year = {2023}
}

@article{dolgov2020approximation,
  title={Approximation and sampling of multivariate probability distributions in the tensor train decomposition},
  author={Dolgov, Sergey and Anaya-Izquierdo, Karim and Fox, Colin and Scheichl, Robert},
  journal={Statistics and Computing},
  volume={30},
  pages={603--625},
  year={2020},
  publisher={Springer}
}

@article{shapeev2016_mtp,
  title={Moment tensor potentials: A class of systematically improvable interatomic potentials},
  author={Shapeev, Alexander V},
  journal={Multiscale Modeling \& Simulation},
  volume={14},
  number={3},
  pages={1153--1173},
  year={2016},
  publisher={SIAM}
}

@article{novikov2020_mlip_2,
  title={The MLIP package: moment tensor potentials with MPI and active learning},
  author={Novikov, Ivan S and Gubaev, Konstantin and Podryabinkin, Evgeny V and Shapeev, Alexander V},
  journal={Machine Learning: Science and Technology},
  volume={2},
  number={2},
  pages={025002},
  year={2020},
  publisher={IOP Publishing}
}

@article{deaven1995_c_clusters,
  title={Molecular geometry optimization with a genetic algorithm},
  author={Deaven, David M and Ho, Kai-Ming},
  journal={Physical review letters},
  volume={75},
  number={2},
  pages={288},
  year={1995},
  publisher={APS}
}

@article{cai2004_brenner,
  title={Structural analysis of carbon clusters by using a global optimization algorithm with Brenner potential},
  author={Cai, Wensheng and Shao, Nan and Shao, Xueguang and Pan, Zhongxiao},
  journal={Journal of Molecular Structure: THEOCHEM},
  volume={678},
  number={1-3},
  pages={113--122},
  year={2004},
  publisher={Elsevier}
}

@inbook{Goreinov2010,
  title = {How to Find a Good Submatrix},
  ISBN = {9789812836021},
  url = {http://dx.doi.org/10.1142/9789812836021_0015},
  DOI = {10.1142/9789812836021_0015},
  booktitle = {Matrix Methods: Theory,  Algorithms and Applications},
  publisher = {WORLD SCIENTIFIC},
  author = {Goreinov,  S. A. and Oseledets,  I. V. and Savostyanov,  D. V. and Tyrtyshnikov,  E. E. and Zamarashkin,  N. L.},
  year = {2010},
  month = apr,
  pages = {247–256}
}

@inbook{CSP_book_intro,

publisher = {John Wiley \& Sons, Ltd},
isbn = {9783527632831},
title = {Introduction: Crystal Structure Prediction, a Formidable Problem},
booktitle = {Modern Methods of Crystal Structure Prediction},
chapter = {},
pages = {11-21},
doiold = {https://doi.org/10.1002/9783527632831.ch},
urlold = {https://onlinelibrary.wiley.com/doi/abs/10.1002/9783527632831.ch},
eprintold = {https://onlinelibrary.wiley.com/doi/pdf/10.1002/9783527632831.ch},
year = {2010},
editor={Artem R. Oganov}
}

@article{USPEX2009,
    author = {Schönborn, Sandro E. and Goedecker, Stefan and Roy, Shantanu and Oganov, Artem R.},
    title = {The performance of minima hopping and evolutionary algorithms for cluster structure prediction},
    journal = {The Journal of Chemical Physics},
    volume = {130},
    number = {14},
    pages = {144108},
    year = {2009},
    doi = {10.1063/1.3097197}
}

@article{Wales1997,
    author = {Wales, David J. and Doye, Jonathan P. K.},
    title = {Global Optimization by Basin-Hopping and the Lowest Energy Structures of Lennard-Jones Clusters Containing up to 110 Atoms},
    journal = {The Journal of Physical Chemistry A},
    volume = {101},
    number = {28},
    pages = {5111-5116},
    year = {1997},
    doi = {10.1021/jp970984n},
    url = {https://doi.org/10.1021/jp970984n}
}

@misc{CSPBench,
      title={CSPBench: a benchmark and critical evaluation of Crystal Structure Prediction}, 
      author={Lai Wei and Sadman Sadeed Omee and Rongzhi Dong and Nihang Fu and Yuqi Song and Edirisuriya M. D. Siriwardane and Meiling Xu and Chris Wolverton and Jianjun Hu},
      year={2024},
      eprint={2407.00733},
      archivePrefix={arXiv},
      primaryClass={cond-mat.mtrl-sci},
      url={https://arxiv.org/abs/2407.00733}, 
}

@book{Wales2003,
    author = {Wales, David J.},
    title = {Energy Landscapes: Applications to Clusters, Biomolecules and Glasses},
    series = {Cambridge Molecular Science},
    publisher = {Cambridge University Press},
    address = {Cambridge},
    year = {2003},
    isbn = {978-0-521-81315-4}
}

@article{USPEX2011,
    author = {Oganov, Artem R. and Lyakhov, Andriy O. and Valle, Mario},
    title = {How evolutionary crystal structure prediction works—and why},
    journal = {Accounts of Chemical Research},
    volume = {44},
    number = {3},
    pages = {227-237},
    year = {2011},
    doi = {10.1021/ar1001318},
    url = {https://doi.org/10.1021/ar1001318}
}

@article{Goedecker2004,
    author = {Goedecker, Stefan},
    title = {Minima hopping: An efficient search method for the global minimum of the potential energy surface of complex molecular systems},
    journal = {The Journal of Chemical Physics},
    volume = {120},
    number = {21},
    pages = {9911-9917},
    year = {2004},
    doi = {10.1063/1.1724816}
}

@inproceedings{
shetty2024generalized,
title={Generalized Policy Iteration using Tensor Approximation for Hybrid Control},
author={Suhan Shetty and Teng Xue and Sylvain Calinon},
booktitle={The Twelfth International Conference on Learning Representations},
year={2024},
url={https://openreview.net/forum?id=csukJcpYDe}
}

@article{shetty2024tt_robotics,
author = {Suhan Shetty and Teguh Lembono and Tobias Löw and Sylvain Calinon},
title ={Tensor train for global optimization problems in robotics},

journal = {The International Journal of Robotics Research},
volume = {43},
number = {6},
pages = {811-839},
year = {2024},
doi = {10.1177/02783649231217527},

URL = { 
    
        https://doi.org/10.1177/02783649231217527
    
    

},
eprint = { 
    
        https://doi.org/10.1177/02783649231217527
    
    

}
}

@book{Bellman:DynamicProgramming,
  author = {Bellman, Richard},
  isbn = {9780486428093},
  publisher = {Dover Publications},
  title = {{Dynamic Programming}},
  year = 1957
}

@article{CPD,
author = {Hitchcock, Frank L.},
title = {The Expression of a Tensor or a Polyadic as a Sum of Products},
journal = {Journal of Mathematics and Physics},
volume = {6},
number = {1-4},
pages = {164-189},
doiold = {https://doi.org/10.1002/sapm192761164},
urlold = {https://onlinelibrary.wiley.com/doi/abs/10.1002/sapm192761164},
eprintold = {https://onlinelibrary.wiley.com/doi/pdf/10.1002/sapm192761164},
year = {1927}
}

@article{Tucker_1966, title={Some Mathematical Notes on Three-Mode Factor Analysis}, volume={31}, DOI={10.1007/BF02289464}, number={3}, journal={Psychometrika}, author={Tucker, Ledyard R}, year={1966}, pages={279–311}}

@article{harshman1970foundations,
  title={Foundations of the {PARAFAC} procedure: Models and conditions for an explanatory multimodal factor analysis},
  author={Harshman, Richard A and others},
  year={1970},
  journal={UCLA Working Papers in Phonetics},
  volume={16},
  pages={1--84},
}

@article{buczyska2015hackbusch,
  title = {The Hackbusch conjecture on tensor formats},
  volume = {104},
  ISSN = {0021-7824},
  DOI = {10.1016/j.matpur.2015.05.002},
  number = {4},
  journal = {Journal de Mathématiques Pures et Appliquées},
  publisher = {Elsevier BV},
  author = {Buczyńska,  Weronika and Buczyński,  Jarosław and Michałek,  Mateusz},
  year = {2015},
  month = oct,
  pages = {749–761}
}

@article{ballani2013black,
  title = {Black box approximation of tensors in hierarchical {T}ucker format},
  volume = {438},
  ISSN = {0024-3795},
  DOI = {10.1016/j.laa.2011.08.010},
  number = {2},
  journal = {Linear Algebra and its Applications},
  publisher = {Elsevier BV},
  author = {Ballani,  Jonas and Grasedyck,  Lars and Kluge,  Melanie},
  year = {2013},
  pages = {639–657}}

@article{TensorRing,
  author       = {Qibin Zhao and
                  Guoxu Zhou and
                  Shengli Xie and
                  Liqing Zhang and
                  Andrzej Cichocki},
  title        = {Tensor Ring Decomposition},
  journal      = {CoRR},
  volume       = {abs/1606.05535},
  year         = {2016},
  url          = {http://arxiv.org/abs/1606.05535},
  eprinttype    = {arXiv},
  eprint       = {1606.05535},
  bibsource    = {dblp computer science bibliography, https://dblp.org}
}

@ARTICLE{Freitas_Gustavo2022-va,
  title    = "{GloMPO} (Globally Managed Parallel Optimization): a tool for
              expensive, black-box optimizations, application to {ReaxFF}
              reparameterizations",
  author   = "Freitas Gustavo, Michael and Verstraelen, Toon",
  journal  = "Journal of Cheminformatics",
  volume   =  14,
  number   =  1,
  pages    = "7",
  month    =  feb,
  year     =  2022
}

@INPROCEEDINGS{TensorRing2,

  author={Zhao, Qibin and Sugiyama, Masashi and Yuan, Longhao and Cichocki, Andrzej},

  booktitle={ICASSP 2019 - 2019 IEEE International Conference on Acoustics, Speech and Signal Processing (ICASSP)}, 

  title={Learning Efficient Tensor Representations with Ring-structured Networks}, 

  year={2019},

  volume={},

  number={},

  pages={8608-8612},

  doi={10.1109/ICASSP.2019.8682231}}

@article{PhysRevA.70.060302,
  title = {Valence-bond states for quantum computation},
  author = {Verstraete, F. and Cirac, J. I.},
  journal = {Phys. Rev. A},
  volume = {70},
  issue = {6},
  pages = {060302},
  numpages = {4},
  year = {2004},
  month = {Dec},
  publisher = {American Physical Society},
  doi = {10.1103/PhysRevA.70.060302},
  url = {https://link.aps.org/doi/10.1103/PhysRevA.70.060302}
}

@misc{verstraete2004,
      title={Renormalization algorithms for Quantum-Many Body Systems in two and higher dimensions}, 
      author={F. Verstraete and J. I. Cirac},
      year={2004},
      eprint={cond-mat/0407066},
      archivePrefix={arXiv},
      primaryClass={cond-mat.str-el},
      url={https://arxiv.org/abs/cond-mat/0407066}, 
}

@article{PhysRevLett.101.110501,
  title = {Class of Quantum Many-Body States That Can Be Efficiently Simulated},
  author = {Vidal, G.},
  journal = {Phys. Rev. Lett.},
  volume = {101},
  issue = {11},
  pages = {110501},
  numpages = {4},
  year = {2008},
  month = {Sep},
  publisher = {American Physical Society},
  doi = {10.1103/PhysRevLett.101.110501},
  url = {https://link.aps.org/doi/10.1103/PhysRevLett.101.110501}
}

@article{PhysRevB.79.144108,
  title = {Algorithms for entanglement renormalization},
  author = {Evenbly, G. and Vidal, G.},
  journal = {Phys. Rev. B},
  volume = {79},
  issue = {14},
  pages = {144108},
  numpages = {20},
  year = {2009},
  month = {Apr},
  publisher = {American Physical Society},
  doi = {10.1103/PhysRevB.79.144108},
  url = {https://link.aps.org/doi/10.1103/PhysRevB.79.144108}
}

@article{PhysRevB.89.235113,
  title = {Scaling of entanglement entropy in the (branching) multiscale entanglement renormalization ansatz},
  author = {Evenbly, G. and Vidal, G.},
  journal = {Phys. Rev. B},
  volume = {89},
  issue = {23},
  pages = {235113},
  numpages = {15},
  year = {2014},
  month = {Jun},
  publisher = {American Physical Society},
  doi = {10.1103/PhysRevB.89.235113},
  url = {https://link.aps.org/doi/10.1103/PhysRevB.89.235113}
}

@article{Batselier2021,
  author    = {Kim Batselier and Andrzej Cichocki and Ngai Wong},
  title     = {MERACLE: Constructive Layer-Wise Conversion of a Tensor Train into a {MERA}},
  journal   = {Communications on Applied Mathematics and Computation},
  year      = {2021},
  volume    = {3},
  number    = {2},
  pages     = {257--279},
  month     = jun,
  doi       = {10.1007/s42967-020-00090-6},
  url       = {https://doi.org/10.1007/s42967-020-00090-6},
  issn      = {2661-8893},
}

@article{cichocki2016tensor,
  title={Tensor networks for dimensionality reduction and large-scale optimization: {Part} 1 low-rank tensor decompositions},
  author={Cichocki, Andrzej and Lee, Namgil and Oseledets, Ivan and Phan, Anh-Huy and Zhao, Qibin and Mandic, Danilo},
  journal={Foundations and Trends in Machine Learning},
  volume={9},
  number={4-5},
  pages={249--429},
  year={2016},
  publisher={Now Publishers Inc. Hanover, MA, USA}
}

@article{cichocki2017tensor,
  title = {Tensor Networks for Dimensionality Reduction and Large-scale Optimization: {Part} 2 Applications and Future Perspectives},
  year = {2017},
  volume = {9},
  journal = {Foundations and Trends in Machine Learning},
  number = {6},
  pages = {431-673},
  author = {Andrzej Cichocki and Anh Phan and Qibin Zhao and Namgil Lee and Ivan Oseledets and Masashi Sugiyama and Danilo Mandic}
}

@article{FIRE,
author = {Rogan, José and Varas, Alejandro and Valdivia, Juan Alejandro and Kiwi, Miguel},
title = {A strategy to find minimal energy nanocluster structures},
journal = {Journal of Computational Chemistry},
volume = {34},
number = {29},
pages = {2548-2556},
doiold = {https://doi.org/10.1002/jcc.23419},
urlold = {https://onlinelibrary.wiley.com/doi/abs/10.1002/jcc.23419},
eprintold = {https://onlinelibrary.wiley.com/doi/pdf/10.1002/jcc.23419},
year = {2013}
}

@article{oseledets2016black,
  title={Black-box solver for multiscale modelling using the QTT format},
  author={Oseledets, Ivan V and Rakhuba, Maxim V and Chertkov, Andrei V},
  journal={Proc. ECCOMAS. Crete Island, Greece},
  year={2016}
}

@book{Nocedal2006Numerical,
  address = {New York},
  author = {Nocedal, Jorge and Wright, Stephen J.},
  biburlold = {https://www.bibsonomy.org/bibtex/255a7423ee5e49aab708f02f70fd05790/gdmcbain},
  comment = {(private-note)`For line search (and optimization-based methods in general), I really like [\~{}]' (Paul T. Bauman, libmesh-users 2012-08-21)},
  doiold = {10.1007/978-0-387-40065-5},
  edition = 2,
  issn = {1431-8598},
  publisher = {Springer},
  series = {Springer Series in Operations Research and Financial Engineering},
  title = {{Numerical Optimization}},
  urlold = {http://dx.doi.org/10.1007/978-0-387-40065-5},
  year = 2006
}

@InProceedings{ECCV2020,
author="Phan, Anh-Huy
and Sobolev, Konstantin
and Sozykin, Konstantin
and Ermilov, Dmitry
and Gusak, Julia
and Tichavsk{\'y}, Petr
and Glukhov, Valeriy
and Oseledets, Ivan
and Cichocki, Andrzej",
editor="Vedaldi, Andrea
and Bischof, Horst
and Brox, Thomas
and Frahm, Jan-Michael",
title="Stable Low-Rank Tensor Decomposition for Compression of Convolutional Neural Network",
booktitle="Computer Vision -- ECCV 2020",
year="2020",
publisher="Springer International Publishing",
address="Cham",
pages="522--539",
isbn="978-3-030-58526-6"
}

@article{Matveev2024,
  author = {Matveev, Segey and Tretyak, Ilya},
  year = {2024},
  month = {11},
  day = {16},
  title = {Nonnegative tensor train for the multicomponent Smoluchowski equation},
  journal = {Computational and Applied Mathematics},
  volume = {44},
  number = {1},
  pages = {35},
  issn = {1807-0302},
  doi = {10.1007/s40314-024-02993-z},
  url = {https://doi.org/10.1007/s40314-024-02993-z}
}

@article{Parsons2005,
  title = {Practical conversion from torsion space to {C}artesian space for \textit{in silico} protein synthesis},
  volume = {26},
  ISSN = {1096-987X},
  url = {http://dx.doi.org/10.1002/jcc.20237},
  DOI = {10.1002/jcc.20237},
  number = {10},
  journal = {Journal of Computational Chemistry},
  publisher = {Wiley},
  author = {Parsons,  Jerod and Holmes,  J. Bradley and Rojas,  J. Maurice and Tsai,  Jerry and Strauss,  Charlie E. M.},
  year = {2005},
  month = may,
  pages = {1063–1068}
}

@article{locatelli2003,
    author = {Locatelli, Marco and Schoen, Fabio},
    title = {Efficient Algorithms for Large Scale Global Optimization: {Lennard-Jones} Clusters},
    journal = {Computational Optimization and Applications},
    volume = {26},
    number = {2},
    pages = {173--190},
    year = {2003},
    month = {11},
    issn = {1573-2894},
    doi = {10.1023/A:1025798414605},
    url = {https://doi.org/10.1023/A:1025798414605}
}

@article{Dolgov+2019+23+38,
url = {https://doi.org/10.1515/cmam-2018-0023},
title = {A Tensor Decomposition Algorithm for Large ODEs with Conservation Laws},
author = {Sergey V. Dolgov},
pages = {23--38},
volume = {19},
number = {1},
journal = {Computational Methods in Applied Mathematics},
doi = {doi:10.1515/cmam-2018-0023},
year = {2019},
lastchecked = {2025-06-18}
}

@misc{dolgov2025tensorcross,
      title={Tensor cross interpolation for global discrete optimization with application to Bayesian network inference}, 
      author={Sergey Dolgov and Dmitry Savostyanov},
      year={2025},
      eprint={2502.12940},
      archivePrefix={arXiv},
      primaryClass={stat.CO},
      url={https://arxiv.org/abs/2502.12940}, 
}

@article{chertkov2023black,
  title={Black box approximation in the tensor train format initialized by {ANOVA} decomposition},
  author={Chertkov, Andrei and Ryzhakov, Gleb and Oseledets, Ivan},
  journal={SIAM Journal on Scientific Computing},
  volume={45},
  number={4},
  pages={A2101--A2118},
  year={2023},
  publisher={SIAM}
}

@article{ballani2012tensor,
  title={Tensor structured evaluation of singular volume integrals},
  author={Ballani, J. and Meszmer, P.},
  journal={Computing and Visualization in Science},
  volume={15},
  number={2},
  pages={75--86},
  year={2012},
  publisher={Springer}
}

@article{litsarev2015fast,
  title={Fast low-rank approximations of multidimensional integrals in ion-atomic collisions modelling},
  author={Litsarev, M. and Oseledets, I.},
  journal={Numerical Linear Algebra with Applications},
  volume={22},
  number={6},
  pages={1147--1160},
  year={2015},
  publisher={Wiley Online Library}
}

@article{vysotsky2021tensor,
  title={Tensor-train numerical integration of multivariate functions with singularities},
  author={Vysotsky, L. and Smirnov, A. and Tyrtyshnikov, E.},
  journal={Lobachevskii Journal of Mathematics},
  volume={42},
  number={7},
  pages={1608--1621},
  year={2021},
  publisher={Springer}
}

@article{khoromskij2010fast,
  title={Fast and accurate tensor approximation of a multivariate convolution with linear scaling in dimension},
  author={Khoromskij, B.},
  journal={Journal of computational and applied mathematics},
  volume={234},
  number={11},
  pages={3122--3139},
  year={2010},
  publisher={Elsevier}
}

@article{rakhuba2015fast,
  title={Fast multidimensional convolution in low-rank tensor formats via cross approximation},
  author={Rakhuba, M. and Oseledets, I.},
  journal={SIAM Journal on Scientific Computing},
  volume={37},
  number={2},
  pages={A565--A582},
  year={2015},
  publisher={SIAM}
}

@article{jin2020ctnn,
  title={{CTNN}: A convolutional tensor-train neural network for multi-task brainprint recognition},
  author={Jin, X. and Tang, J. and Kong, X. and Peng, Y. and Cao, J. and Zhao, Q. and Kong, W.},
  journal={IEEE Transactions on Neural Systems and Rehabilitation Engineering},
  volume={29},
  pages={103--112},
  year={2020},
  publisher={IEEE}
}

@article{fernandez2022learning,
  title={Learning {F}eynman diagrams with tensor trains},
  author={Fern{\'a}ndez, Y. and Jeannin, M. and Dumitrescu, P. and Kloss, T. and Kaye, J. and Parcollet, O. and Waintal, X.},
  journal={Physical Review X},
  volume={12},
  number={4},
  pages={041018},
  year={2022},
  publisher={APS}
}

@article{erpenbeck2023tensor,
  title={A tensor train continuous time solver for quantum impurity models},
  author={Erpenbeck, A. and Lin, W. and Blommel, T. and Zhang, L. and Iskakov, S. and Bernheimer, L. and N{\'u}{\~n}ez-Fern{\'a}ndez, Y. and Cohen, G. and Parcollet, O. and Waintal, X. and others},
  journal={arXiv preprint arXiv:2303.11199},
  year={2023}
}

@article{shinaoka2023multiscale,
  title={Multiscale space-time ansatz for correlation functions of quantum systems based on quantics tensor trains},
  author={Shinaoka, H. and Wallerberger, M. and Murakami, Y. and Nogaki, K. and Sakurai, R. and Werner, P. and Kauch, A.},
  journal={Physical Review X},
  volume={13},
  number={2},
  pages={021015},
  year={2023},
  publisher={APS}
}

@article{glau2020low,
  title={Low-rank tensor approximation for {Chebyshev} interpolation in parametric option pricing},
  author={Glau, K. and Kressner, D. and Statti, F.},
  journal={SIAM Journal on Financial Mathematics},
  volume={11},
  number={3},
  pages={897--927},
  year={2020},
  publisher={SIAM}
}

@inproceedings{richter2021solving,
  title={Solving high-dimensional parabolic {PDE}s using the tensor train format},
  author={Richter, L. and Sallandt, L. and N{\"u}sken, N.},
  booktitle={International Conference on Machine Learning},
  pages={8998--9009},
  year={2021},
  organization={PMLR}
}

@article{bayer2023pricing,
  title={Pricing high-dimensional {B}ermudan options with hierarchical tensor formats},
  author={Bayer, C. and Eigel, M. and Sallandt, L. and Trunschke, P.},
  journal={SIAM Journal on Financial Mathematics},
  volume={14},
  number={2},
  pages={383--406},
  year={2023},
  publisher={SIAM}
}

@article{ahmadiasl2021cross,
  title={Cross tensor approximation methods for compression and dimensionality reduction},
  author={Ahmadi-Asl, S. and Caiafa, C. and Cichocki, A. and Phan, A. and Tanaka, T. and Oseledets, I. and Wang, J.},
  journal={IEEE Access},
  volume={9},
  pages={150809--150838},
  year={2021},
  publisher={IEEE}
}

@article{lee2021qttnet,
  title={{QTTNet}: Quantized tensor train neural networks for 3D object and video recognition},
  author={Lee, D. and Wang, D. and Yang, Y. and Deng, L. and Zhao, G. and Li, G.},
  journal={Neural Networks},
  volume={141},
  pages={420--432},
  year={2021},
  publisher={Elsevier}
}

@article{qiu2022efficient,
  title={Efficient tensor robust {PCA} under hybrid model of tucker and tensor train},
  author={Qiu, Y. and Zhou, G. and Huang, Z. and Zhao, Q. and Xie, S.},
  journal={IEEE Signal Processing Letters},
  volume={29},
  pages={627--631},
  year={2022},
  publisher={IEEE}
}

@inproceedings{tjandra2017compressing,
  title={Compressing recurrent neural network with tensor train},
  author={Tjandra, A. and Sakti, S. and Nakamura, S.},
  booktitle={2017 International Joint Conference on Neural Networks (IJCNN)},
  pages={4451--4458},
  year={2017},
  organization={IEEE}
}

@article{novikov2020tensor,
  title={Tensor train decomposition on tensorflow (t3f)},
  author={Novikov, A. and Izmailov, P. and Khrulkov, V. and Figurnov, M. and Oseledets, I.},
  journal={The Journal of Machine Learning Research},
  volume={21},
  number={1},
  pages={1105--1111},
  year={2020},
  publisher={JMLRORG}
}

@article{wang2021nonlinear,
  title={Nonlinear tensor train format for deep neural network compression},
  author={Wang, D. and Zhao, G. and Chen, H. and Liu, Z. and Deng, L. and Li, G.},
  journal={Neural Networks},
  volume={144},
  pages={320--333},
  year={2021},
  publisher={Elsevier}
}

@inproceedings{chen2019support,
  title={A support tensor train machine},
  author={Chen, C. and Batselier, K. and Ko, C. and Wong, N.},
  booktitle={2019 International Joint Conference on Neural Networks (IJCNN)},
  pages={1--8},
  year={2019},
  organization={IEEE}
}

@article{kour2023efficient,
  title={Efficient structure-preserving support tensor train machine},
  author={Kour, K. and Dolgov, S. and Stoll, M. and Benner, P.},
  journal={Journal of Machine Learning Research},
  volume={24},
  number={4},
  pages={1--22},
  year={2023}
}

@article{liu2023tensor,
  title={Tensor networks for unsupervised machine learning},
  author={Liu, J. and Li, S. and Zhang, J. and Zhang, P.},
  journal={Physical Review E},
  volume={107},
  number={1},
  pages={L012103},
  year={2023},
  publisher={APS}
}

@article{chertkov2016robust,
  title         = {Robust discretization in quantized tensor train format for elliptic problems in two dimensions},
  author        = {Chertkov, A. and Oseledets, I. and Rakhuba, M.},
  journal       = {arXiv preprint arXiv:1612.01166},
  year          = {2016}
}

@article{chertkov2021solution,
    author      = {Chertkov, A. and Oseledets, I.},
    title       = {Solution of the {F}okker–{P}lanck equation by cross approximation method in the tensor train format},
    journal     = {Frontiers in Artificial Intelligence},
    volume      = {4},
    year        = {2021}
}

@article{chertkov2024translate,
  title={Translate your gibberish: black-box adversarial attack on machine translation systems},
  author={Chertkov, Andrei and Tsymboi, Olga and Pautov, Mikhail and Oseledets, Ivan},
  journal={Journal of Mathematical Sciences},
  volume={285},
  number={2},
  pages={221--233},
  year={2024},
  publisher={Springer}
}

@article{chertkov2023tensor,
  title={Tensor Train Decomposition for Adversarial Attacks on Computer Vision Models},
  author={Chertkov, Andrei and Oseledets, Ivan},
  journal={arXiv preprint arXiv:2312.12556},
  year={2023}
}

@article{sozykin2025high,
  title={High-dimensional Optimization with Low Rank Tensor Sampling and Local Search},
  author={Sozykin, Konstantin and Chertkov, Andrei and Phan, Anh-Huy and Oseledets, Ivan and Ryzhakov, Gleb},
  journal={arXiv preprint arXiv:2505.12383},
  year={2025}
}

@article{pospelov2025fast,
  title={Fast gradient-free activation maximization for neurons in spiking neural networks},
  author={Pospelov, Nikita and Chertkov, Andrei and Beketov, Maxim and Oseledets, Ivan and Anokhin, Konstantin},
  journal={Neurocomputing},
  volume={618},
  pages={129070},
  year={2025},
  publisher={Elsevier}
}

@inproceedings{lee2024tt,
  title={TT-SNN: tensor train decomposition for efficient spiking neural network training},
  author={Lee, Donghyun and Yin, Ruokai and Kim, Youngeun and Moitra, Abhishek and Li, Yuhang and Panda, Priyadarshini},
  booktitle={2024 Design, Automation \& Test in Europe Conference \& Exhibition (DATE)},
  pages={1--6},
  year={2024},
  organization={IEEE}
}

@article{chertkov2023tensorExtrema,
  title={Tensor Extrema Estimation Via Sampling: A New Approach for Determining Minimum/Maximum Elements},
  author={Chertkov, Andrei and Ryzhakov, Gleb and Novikov, Georgii and Oseledets, Ivan},
  journal={Computing in Science \& Engineering},
  volume={25},
  number={5},
  pages={14--25},
  year={2023},
  publisher={IEEE}
}

@article{ryzhakov2024black,
  title={Black-Box Approximation and Optimization with Hierarchical Tucker Decomposition},
  author={Ryzhakov, Gleb and Chertkov, Andrei and Basharin, Artem and Oseledets, Ivan},
  journal={arXiv preprint arXiv:2402.02890},
  year={2024}
}

@article{mikhalev2018rectangular,
  title={Rectangular maximum-volume submatrices and their applications},
  author={Mikhalev, Aleksandr and Oseledets, Ivan V},
  journal={Linear Algebra and its Applications},
  volume={538},
  pages={187--211},
  year={2018},
  publisher={Elsevier}
}

@article{Orus2019,
  author    = {Roman Orus},
  title     = {Tensor networks for complex quantum systems},
  journal   = {Nature Reviews Physics},
  year      = {2019},
  volume    = {1},
  number    = {9},
  pages     = {538--550},
  doi       = {10.1038/s42254-019-0086-7},
  url       = {https://doi.org/10.1038/s42254-019-0086-7},
  issn      = {2522-5820}
}

@misc{xue2025unifyingrobotoptimizationmonte,
      title={Unifying Robot Optimization: Monte Carlo Tree Search with Tensor Factorization}, 
      author={Teng Xue and Amirreza Razmjoo and Yan Zhang and Sylvain Calinon},
      year={2025},
      eprint={2507.04949},
      archivePrefix={arXiv},
      primaryClass={cs.RO},
      url={https://arxiv.org/abs/2507.04949}, 
}

@article{GAS,
    author = {Chen, Zhanghui and Jiang, Xiangwei and Li, Jingbo and Li, Shushen},
    title = {A sphere-cut-splice crossover for the evolution of cluster structures},
    journal = {The Journal of Chemical Physics},
    volume = {138},
    number = {21},
    pages = {214303},
    year = {2013},
    month = {06},
    issn = {0021-9606},
    doi = {10.1063/1.4807091},
    url = {https://doi.org/10.1063/1.4807091},
    eprint = {},
}

@article{WalesLJ38,
    author = {Doye, Jonathan P. K. and Miller, Mark A. and Wales, David J.},
    title = {The double-funnel energy landscape of the 38-atom Lennard-Jones cluster},
    journal = {The Journal of Chemical Physics},
    volume = {110},
    number = {14},
    pages = {6896-6906},
    year = {1999},
    month = {04},
    issn = {0021-9606},
    doi = {10.1063/1.478595},
    url = {https://doi.org/10.1063/1.478595},
}

@incollection{zamarashkin2010-maxvol,
  title={How to find a good submatrix},
  author={Goreinov, Sergei A and Oseledets, Ivan V and Savostyanov, Dimitry V and Tyrtyshnikov, Eugene E and Zamarashkin, Nikolay L},
  booktitle={Matrix Methods: Theory, Algorithms And Applications: Dedicated to the Memory of Gene Golub},
  pages={247--256},
  year={2010},
  publisher={World Scientific}
}

@article{Podryabinkin2017,
title = {Active learning of linearly parametrized interatomic potentials},
journal = {Computational Materials Science},
volume = {140},
pages = {171-180},
year = {2017},
issn = {0927-0256},
doi = {https://doi.org/10.1016/j.commatsci.2017.08.031},
url = {https://www.sciencedirect.com/science/article/pii/S0927025617304536},
author = {Evgeny V. Podryabinkin and Alexander V. Shapeev},
}

@book{Amari2000Methods,
  author    = {Amari, Shun-ichi and Nagaoka, Hiroshi},
  title     = {Methods of Information Geometry},
  series    = {Translations of Mathematical Monographs},
  volume    = {191},
  publisher = {American Mathematical Society},
  year      = {2000},
  isbn      = {0821843028},
  pages     = {206},
  note      = {Translated by Daishi Harada},
  url       = {https://bookstore.ams.org/mmono-191}
}

@Article{e22101100,
AUTHOR = {Nielsen, Frank},
TITLE = {An Elementary Introduction to Information Geometry},
JOURNAL = {Entropy},
VOLUME = {22},
YEAR = {2020},
NUMBER = {10},
ARTICLE-NUMBER = {1100},
URL = {https://www.mdpi.com/1099-4300/22/10/1100},
PubMedID = {33286868},
ISSN = {1099-4300},
DOI = {10.3390/e22101100}
}
\appendix
\section{Supplementary materials}
\label{sec:supp}

This appendix presents the complete numerical results underlying the performance analysis in the main text. 
It includes detailed tables for the PROTES optimizer with all tested encoding schemes (Simple Relative, Constrained Relative, and Direct Cartesian) across Lennard-Jones clusters from 5 to 45 atoms, 
as well as results for the TTOpt baseline. 
The data enable a granular examination of computational scaling, success rates, and the interaction between encoding choice and initialization strategy, 
providing the empirical foundation for the efficiency-reliability trade-off discussed in Section~\ref{sec:results_lj}.
Additionally, we detail the hyperparameter settings used for both PROTES and TTOpt optimizers.

\subsection{Hyperparameter Selection}

\begin{figure}[h]
    \centering
    \includegraphics[width=0.7\linewidth]{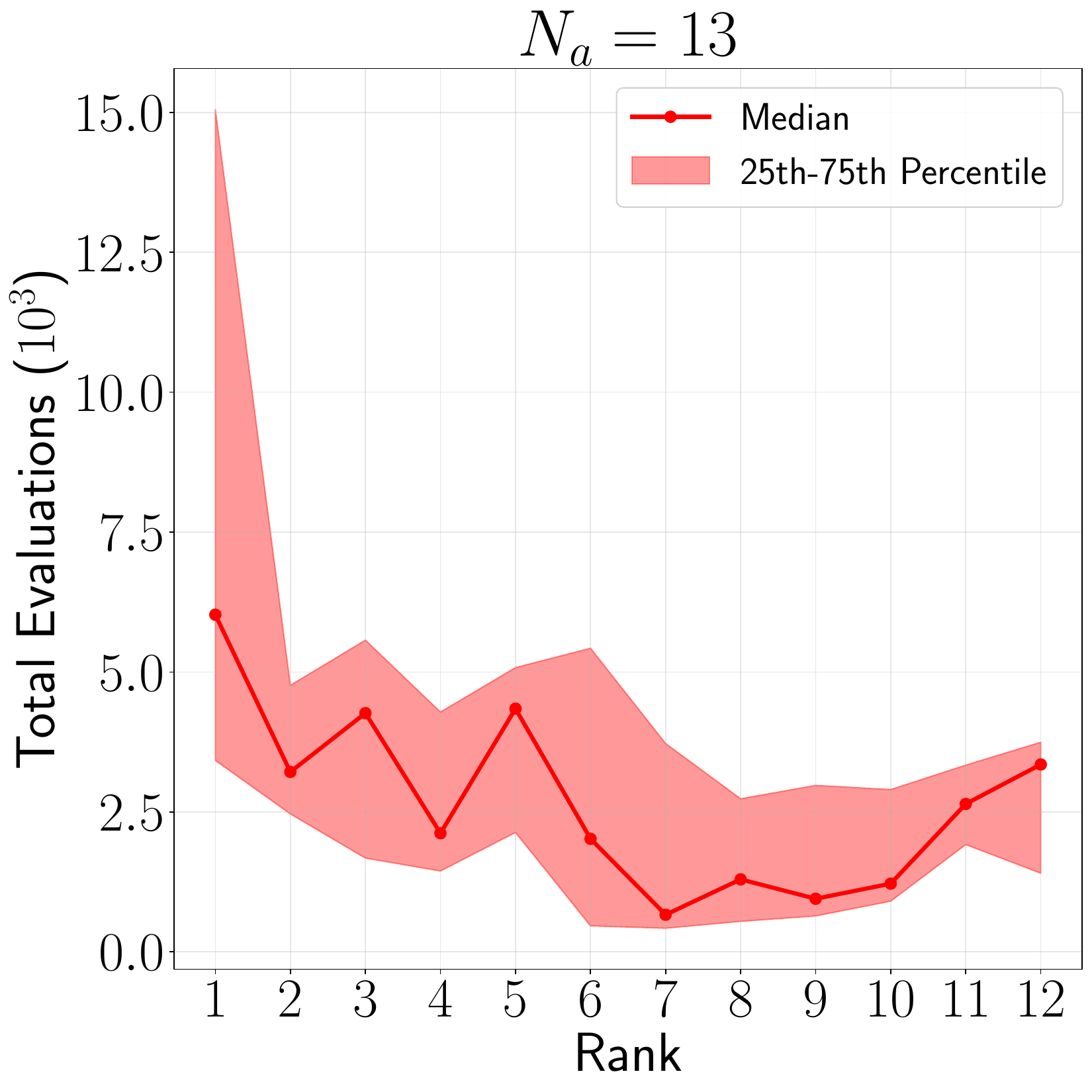}
    \caption{
        TT-rank sensitivity analysis for LJ$_{13}$ indicates optimal performance at $R = 7$.
    }
    \label{fig:rank_sensitivity}
\end{figure}

All experiments were conducted with the following consistent parameters:
\begin{itemize}
    \item \texttt{Potential Parameters}: $\epsilon = 1$, $\sigma = 1$ (reduced units)
    \item \texttt{Search Domain}: $\Omega = [-2, 2]^{3M}$ for direct encoding; appropriate physical bounds for relative encodings
    \item \texttt{Grid Resolution}: $N_i = 16-32$ points per dimension depending on encoding scheme
    \item \texttt{TT-Ranks}: The TT-rank for all methods is empirically selected to $R = 7$ based on sensitivity analysis of a LJ$_{13}$(Figure~\ref{fig:rank_sensitivity})
    \item \texttt{Local Optimizer}: L-BFGS-B with potential forces as the analytical gradients
    \item \texttt{Convergence Tolerance}: $10^{-4}$.
\end{itemize}

\subsection{Statistical analysis of TTOpt}

Table~\ref{tab:ttopt_results} provides a detailed statistical analysis of the \texttt{TTOpt algorithm} using direct encoding, averaged over 30 independent runs.
The results demonstrate exceptional statistical robustness, with all runs achieving 100\% success rates across the tested clusters (5 to 26 atoms) and energy standard deviations in the range of $10^{-15}$ to $10^{-14}$ eV.
The computational effort exhibits predictable polynomial scaling, although the absolute costs exceed those of PROTES with advanced encodings.
The increased variance in function calls for specific clusters (e.g., LJ$_{16}$, LJ$_{19}$) suggests that these systems possess particularly challenging optimization landscapes, resulting in greater variability in the number of required iterations.

\begin{table}[t!]
    \centering
    \small
    \caption{TTOpt performance with direct encoding for Lennard-Jones clusters. Results are averaged across 30 random seeds, all achieving 100\% success rate (SR$_t$). The number of calls represents total evaluations TC (TT function calls including local optimization). RE denotes the relative error in energy.}
    \label{tab:ttopt_results}
    \begin{tabular}{|c|r|r|r|}
    \hline
    $N_a$ & \multicolumn{1}{c|}{TC} & \multicolumn{1}{c|}{RE} & \multicolumn{1}{c|}{SR$_t$ (\%)} \\
    \hline
    5  & 3.16E+02 & 2.53E-12 & 100 \\
    7  & 9.80E+03 & 9.69E-12 & 100 \\
    8  & 8.23E+03 & 1.99E-10 & 100 \\
    9  & 2.95E+03 & 5.40E-11 & 100 \\
    10 & 9.18E+03 & 1.12E-10 & 100 \\
    11 & 4.59E+03 & 3.67E-10 & 100 \\
    13 & 1.07E+04 & 1.57E-10 & 100 \\
    15 & 1.83E+04 & 1.70E-10 & 100 \\
    16 & 1.47E+04 & 6.81E-03 & 100 \\
    17 & 2.16E+04 & 5.01E-10 & 100 \\
    19 & 2.22E+04 & 3.79E-10 & 100 \\
    20 & 4.60E+04 & 6.51E-10 & 100 \\
    24 & 3.12E+04 & 6.09E-10 & 100 \\
    25 & 9.28E+04 & 6.92E-09 & 100 \\
    26 & 1.95E+04 & 4.63E-10 & 100 \\
    \hline
    \end{tabular}
\end{table}

\subsection{Detailed performance tables}

The supplementary tables present complete numerical results for the PROTES optimizer, organized by encoding scheme and cluster type. Rather than testing every consecutive size, our experiments focus on representative clusters spanning different structural classes and optimization difficulty regimes within the 5–45 atom range. This includes magic-number clusters with symmetric, funnel-like landscapes, as well as challenging non-magic numbers with frustrated or multi-funnel energy surfaces. The comprehensive data provided underpins the strategic analysis of the efficiency-reliability trade-off presented in Section~\ref{sec:results_lj}.

\subsubsection{PROTES with relative encodings}

Table~\ref{table:protes_simple_relative} summarizes the results for the Simple Relative encoding (relative encoding with constant distances).
The data reveal a consistent polynomial scaling of computational cost with respect to system size, where PROTES calls (PC) typically range from $10^1$--$10^2$ for small clusters ($N_a \leq 15$) to $10^4$--$10^5$ for larger systems ($N_a \geq 30$).
Success rates persist at 100\% for clusters up to $N_a = 28$, after which they decline gradually, falling to 0\% for the challenging LJ$_{45}$ system despite significant computational effort.

Table~\ref{table:protes_constrained_relative} presents the performance metrics for the Constrained Relative encoding (relative encoding with angular constraints).
While this encoding exhibits similar scaling trends, it generally incurs higher computational costs than the Simple Relative encoding for medium-sized clusters.
Nevertheless, it demonstrates superior robustness for the largest systems, maintaining a 33\% success rate for LJ$_{45}$, whereas the Simple Relative encoding fails completely.
This observation reinforces the conclusion that more constrained encodings enhance reliability at the expense of computational efficiency.

\begin{table}[H]
 \tiny
  \caption{PROTES Relative encoding with constant distances. $N_{a}$: number of atoms; PC: PROTES calls (energy evaluations); Local calls total (LCT): total local optimization evaluations; Local calls last (LCL): last local optimization evaluations; RE: relative error; $SR_t$: success rate (\%).}
  \label{table:protes_simple_relative}
  \centering
  \begin{tabular}{|c|r|r|r|r|r|}
    \hline
    $N_{a}$ & PC & LCT & LCL & RE & $SR_t$ \\
    \hline
5 & 2.02E+01 & 6.73E+02 & 1.39E+02 & 2.62E-12 & 100 \\
6 & 3.68E+02 & 5.31E+03 & 1.39E+02 & 1.65E-12 & 100 \\
7 & 4.63E+01 & 2.64E+03 & 2.47E+02 & 1.07E-11 & 100 \\
8 & 1.07E+02 & 1.34E+03 & 1.37E+02 & 2.19E-10 & 100 \\
9 & 4.10E+01 & 3.52E+03 & 1.88E+02 & 1.27E-10 & 100 \\
10 & 2.08E+02 & 2.90E+03 & 1.19E+02 & 3.35E-12 & 100 \\
11 & 2.06E+02 & 1.13E+04 & 2.86E+02 & 5.59E-11 & 100 \\
12 & 2.22E+02 & 4.33E+03 & 1.61E+02 & 9.71E-11 & 100 \\
13 & 1.04E+02 & 5.44E+03 & 3.89E+02 & 4.16E-10 & 100 \\
14 & 1.78E+02 & 3.83E+03 & 1.69E+02 & 1.17E-10 & 100 \\
15 & 5.60E+01 & 3.93E+03 & 3.45E+02 & 3.46E-10 & 100 \\
16 & 1.63E+02 & 9.65E+03 & 3.48E+02 & 1.30E-10 & 100 \\
17 & 1.04E+03 & 3.27E+04 & 1.54E+02 & 1.87E-10 & 100 \\
18 & 6.66E+03 & 1.50E+05 & 1.59E+02 & 1.76E-10 & 100 \\
19 & 6.31E+02 & 1.38E+04 & 1.19E+02 & 2.80E-10 & 100 \\
20 & 3.71E+02 & 1.28E+04 & 2.88E+02 & 7.37E-10 & 100 \\
21 & 6.22E+03 & 1.12E+05 & 2.14E+02 & 6.73E-10 & 100 \\
22 & 2.78E+03 & 6.44E+04 & 2.03E+02 & 4.77E-10 & 100 \\
23 & 7.67E+03 & 1.14E+05 & 2.45E+02 & 1.30E-03 & 91 \\
24 & 1.39E+04 & 1.90E+05 & 4.07E+02 & 4.47E-10 & 100 \\
25 & 2.42E+02 & 1.16E+04 & 2.57E+02 & 6.94E-10 & 100 \\
26 & 2.15E+03 & 2.32E+05 & 5.08E+02 & 1.03E-09 & 100 \\
27 & 4.78E+03 & 2.39E+05 & 2.72E+02 & 8.77E-10 & 100 \\
28 & 2.42E+03 & 7.52E+04 & 3.93E+02 & 3.68E-10 & 100 \\
29 & 2.44E+04 & 8.22E+05 & 2.72E+02 & 4.77E-03 & 33 \\
30 & 2.18E+04 & 6.01E+05 & 2.87E+02 & 9.96E-04 & 25 \\
31 & 2.59E+04 & 4.69E+05 & 2.81E+02 & 3.02E-03 & 16 \\
32 & 2.85E+04 & 1.48E+06 & 4.11E+02 & 7.93E-03 & 27 \\
33 & 2.21E+04 & 9.46E+05 & 3.68E+02 & 2.49E-03 & 41 \\
34 & 2.58E+04 & 1.05E+06 & 3.64E+02 & 6.41E-03 & 0 \\
38 & 4.89E+04 & 3.78E+06 & 4.46E+02 & 8.27E-03 & 16 \\
40 & 5.28E+04 & 2.86E+06 & 5.34E+02 & 3.37E-03 & 66 \\
41 & 7.13E+04 & 4.95E+06 & 6.89E+02 & 9.29E-03 & 0 \\
43 & 4.32E+04 & 5.04E+06 & 6.50E+02 & 5.15E-03 & 50 \\
45 & 5.56E+04 & 3.34E+06 & 7.67E+02 & 6.41E-03 & 0 \\
    \hline
  \end{tabular}
\end{table}

\begin{table}[h]
 \tiny
  \caption{PROTES Relative encoding with constant distances and restrictions on angles. $N_{a}$: number of atoms; PC: PROTES calls (energy evaluations); Local calls total (LCT): total local optimization evaluations; Local calls last (LCL): last local optimization evaluations; RE: relative error; $SR_t$: success rate (\%).}
  \label{table:protes_constrained_relative}
  \centering
  \begin{tabular}{|c|r|r|r|r|r|}
    \hline
    $N_{a}$ & PC & LCT & LCL & RE & $SR_t$ \\
    \hline
5 & 3.82E+01 & 6.48E+02 & 1.02E+02 & 2.60E-12 & 100 \\
6 & 9.88E+01 & 2.25E+03 & 1.64E+02 & 8.07E-13 & 100 \\
7 & 1.85E+03 & 1.81E+03 & 1.23E+02 & 6.12E-11 & 100 \\
8 & 8.20E+01 & 2.38E+03 & 6.52E+01 & 2.19E-10 & 100 \\
9 & 5.00E+02 & 9.58E+03 & 1.45E+02 & 3.12E-10 & 100 \\
10 & 1.72E+02 & 4.15E+03 & 8.44E+01 & 1.97E-12 & 100 \\
11 & 4.14E+03 & 6.60E+04 & 2.27E+02 & 1.19E-10 & 100 \\
12 & 1.27E+02 & 3.28E+03 & 1.73E+02 & 9.64E-11 & 100 \\
13 & 2.54E+04 & 7.16E+05 & 4.07E+02 & 3.91E-10 & 100 \\
14 & 9.48E+01 & 3.75E+03 & 1.48E+02 & 1.17E-10 & 100 \\
15 & 1.49E+04 & 5.05E+03 & 2.69E+02 & 2.34E-10 & 100 \\
16 & 1.25E+03 & 1.23E+04 & 3.04E+02 & 2.88E-10 & 100 \\
17 & 4.16E+02 & 1.83E+04 & 1.18E+02 & 1.88E-10 & 100 \\
18 & 7.13E+03 & 2.62E+05 & 1.60E+02 & 1.76E-10 & 100 \\
19 & 3.23E+02 & 1.63E+04 & 1.53E+02 & 2.81E-10 & 100 \\
20 & 1.27E+04 & 1.67E+04 & 2.85E+02 & 7.01E-10 & 100 \\
21 & 4.63E+02 & 2.77E+04 & 1.88E+02 & 6.72E-10 & 100 \\
22 & 2.22E+03 & 1.01E+05 & 2.30E+02 & 4.77E-10 & 100 \\
23 & 1.95E+03 & 8.57E+04 & 2.02E+02 & 4.48E-10 & 100 \\
24 & 5.25E+04 & 7.43E+05 & 5.49E+02 & 6.28E-10 & 100 \\
25 & 2.18E+02 & 1.23E+04 & 2.69E+02 & 6.93E-10 & 100 \\
26 & 1.75E+04 & 1.02E+06 & 2.63E+02 & 1.03E-09 & 100 \\
27 & 5.28E+03 & 3.53E+05 & 2.44E+02 & 8.76E-10 & 100 \\
28 & 5.99E+03 & 2.92E+05 & 3.20E+02 & 3.68E-10 & 100 \\
29 & 2.48E+04 & 1.47E+06 & 2.76E+02 & 5.17E-03 & 25 \\
30 & 1.82E+04 & 8.49E+05 & 2.62E+02 & 4.64E-04 & 50 \\
31 & 2.33E+04 & 1.73E+06 & 3.26E+02 & 1.96E-03 & 16 \\
32 & 1.43E+04 & 7.88E+05 & 3.68E+02 & 6.31E-03 & 8 \\
33 & 2.34E+04 & 2.25E+06 & 3.54E+02 & 4.05E-03 & 25 \\
34 & 3.56E+04 & 3.57E+06 & 3.84E+02 & 1.53E-03 & 33 \\
38 & 5.45E+04 & 6.62E+06 & 3.85E+02 & 4.20E-03 & 33 \\
40 & 3.63E+04 & 4.88E+06 & 6.07E+02 & 1.49E-03 & 66 \\
41 & 3.31E+04 & 5.46E+06 & 4.64E+02 & 2.82E-03 & 0 \\
43 & 2.91E+04 & 6.06E+06 & 6.64E+02 & 5.03E-03 & 36 \\
45 & 5.23E+04 & 1.11E+07 & 1.61E+03 & 4.44E-03 & 33 \\
    \hline
  \end{tabular}
\end{table}

\subsubsection{Impact of the initialization strategy}

Tables~\ref{table:protes_constrained_relative_anostic} and \ref{table:protes_constrained_relative_agnostic} facilitate a direct comparison of initialization strategies for both relative encoding schemes.
The data clearly illustrate a size-dependent transition in the optimal strategy: agnostic initialization proves superior for $N_a \leq 26$ due to significantly lower computational costs, whereas physically-constrained initialization becomes essential for $N_a \geq 30$ to achieve successful convergence.
For intermediate clusters (e.g., LJ$_{18}$--LJ$_{25}$), the choice represents a distinct trade-off, requiring the user to prioritize either efficiency (agnostic) or reliability (constrained).

\begin{table}[H]
 \tiny
  \caption{PROTES Relative encoding with constant distances and without physically-constrained  initialization. $N_{a}$: number of atoms; PC: PROTES +  calls (energy evaluations); Local calls total (LCT): total local optimization evaluations; Local calls last (LCL): last local optimization evaluations; RE: relative error; $SR_t$: success rate (\%).}
  \label{table:protes_constrained_relative_anostic}
  \centering
  \begin{tabular}{|c|r|r|r|r|r|}
    \hline
    $N_{a}$ & PC & LCT & LCL & RE & $SR_t$ \\
    \hline
5 & 4.83E+01 & 1.35E+03 & 3.25E+02 & 2.48E-12 & 100 \\
6 & 2.93E+02 & 4.62E+03 & 1.49E+02 & 1.01E-12 & 100 \\
7 & 1.72E+01 & 5.31E+02 & 1.07E+02 & 2.94E-11 & 100 \\
8 & 3.10E+01 & 5.90E+02 & 1.20E+02 & 2.19E-10 & 100 \\
9 & 3.62E+01 & 9.30E+02 & 1.08E+02 & 2.23E-11 & 100 \\
10 & 2.82E+02 & 4.29E+03 & 1.37E+02 & 3.65E-12 & 100 \\
11 & 2.52E+02 & 4.37E+03 & 1.04E+02 & 1.61E-10 & 100 \\
12 & 4.74E+02 & 7.46E+03 & 1.64E+02 & 9.65E-11 & 100 \\
13 & 1.26E+02 & 3.08E+03 & 1.34E+02 & 1.71E-11 & 100 \\
14 & 1.97E+02 & 4.04E+03 & 1.71E+02 & 1.17E-10 & 100 \\
15 & 9.42E+01 & 2.11E+03 & 1.61E+02 & 2.99E-10 & 100 \\
16 & 3.61E+02 & 8.54E+03 & 1.32E+02 & 2.35E-10 & 100 \\
17 & 1.13E+03 & 1.92E+04 & 1.47E+02 & 1.88E-10 & 100 \\
18 & 1.24E+04 & 2.38E+05 & 1.71E+02 & 1.77E-10 & 100 \\
19 & 6.22E+02 & 1.15E+04 & 1.49E+02 & 2.82E-10 & 100 \\
20 & 1.07E+03 & 1.76E+04 & 1.81E+02 & 1.04E-09 & 100 \\
21 & 2.45E+03 & 3.81E+04 & 2.06E+02 & 6.72E-10 & 100 \\
22 & 8.28E+03 & 9.50E+04 & 2.05E+02 & 4.76E-10 & 100 \\
23 & 1.02E+04 & 1.05E+05 & 1.92E+02 & 4.47E-10 & 100 \\
24 & 8.14E+03 & 7.43E+04 & 2.29E+02 & 6.95E-10 & 100 \\
25 & 1.52E+03 & 3.97E+04 & 2.92E+02 & 6.94E-10 & 100 \\
26 & 1.82E+04 & 3.90E+05 & 2.30E+02 & 3.69E-03 & 66 \\
27 & 1.17E+04 & 2.81E+05 & 2.67E+02 & 8.76E-10 & 100 \\
    \hline
  \end{tabular}
\end{table}

\begin{table}[H]
 \tiny
  \caption{PROTES Relative encoding with constant distances and restrictions on angles and without physically-constrained  initialization. $N_{a}$: number of atoms; PC: PROTES +  calls (energy evaluations); Local calls total (LCT): total local optimization evaluations; Local calls last (LCL): last local optimization evaluations; RE: relative error; $SR_t$: success rate (\%).}
  \label{table:protes_constrained_relative_agnostic}
  \centering
  \begin{tabular}{|c|r|r|r|r|r|}
    \hline
    $N_{a}$ & PC & LCT & LCL & RE & SR (\%) \\
    \hline
5 & 5.26E+01 & 3.24E+03 & 1.51E+02 & 2.70E-12 & 100 \\
6 & 8.58E+01 & 2.16E+03 & 1.60E+02 & 9.54E-13 & 100 \\
7 & 5.40E+01 & 1.94E+03 & 1.47E+02 & 3.01E-11 & 100 \\
8 & 3.11E+01 & 1.08E+03 & 8.40E+01 & 2.19E-10 & 100 \\
9 & 4.86E+01 & 1.23E+03 & 1.35E+02 & 2.31E-11 & 100 \\
10 & 2.12E+02 & 5.24E+03 & 8.68E+01 & 2.36E-12 & 100 \\
11 & 3.01E+02 & 6.59E+03 & 1.11E+02 & 1.62E-10 & 100 \\
12 & 7.32E+01 & 2.07E+03 & 1.91E+02 & 3.79E-03 & 91 \\
13 & 5.61E+01 & 1.64E+03 & 1.22E+02 & 1.67E-11 & 100 \\
14 & 6.27E+01 & 2.54E+03 & 1.51E+02 & 1.17E-10 & 100 \\
15 & 5.11E+01 & 2.10E+03 & 1.53E+02 & 2.99E-10 & 100 \\
16 & 1.88E+02 & 6.32E+03 & 1.66E+02 & 2.34E-10 & 100 \\
17 & 4.36E+02 & 1.34E+04 & 1.50E+02 & 1.87E-10 & 100 \\
18 & 1.67E+02 & 6.42E+03 & 1.64E+02 & 3.70E-03 & 0 \\
19 & 2.74E+02 & 9.57E+03 & 1.42E+02 & 1.27E-02 & 50 \\
20 & 2.62E+02 & 8.75E+03 & 2.28E+02 & 3.13E-03 & 75 \\
21 & 2.16E+03 & 7.25E+04 & 1.83E+02 & 6.72E-10 & 100 \\
22 & 7.80E+01 & 3.26E+03 & 1.46E+02 & 4.12E-03 & 50 \\
23 & 1.66E+03 & 4.19E+04 & 1.96E+02 & 3.12E-03 & 80 \\
24 & 6.51E+03 & 1.33E+05 & 2.35E+02 & 5.25E-05 & 85 \\
25 & 2.26E+02 & 1.06E+04 & 1.78E+02 & 2.79E-03 & 50 \\
26 & 1.63E+04 & 9.07E+05 & 2.12E+02 & 8.58E-03 & 25 \\
27 & 5.95E+03 & 2.47E+05 & 3.03E+02 & 5.69E-04 & 90 \\
    \hline
  \end{tabular}
\end{table}

\subsubsection{Baseline performance with direct encoding}

Table~\ref{table:protes_naive_encoding} documents the performance of PROTES utilizing \texttt{Direct Cartesian encoding}, serving as a baseline that emphasizes the critical importance of physically-informed discretization.
The results indicate substantially higher computational costs and rapidly declining success rates beyond $N_a = 10$, rendering the method computationally prohibitive for $N_a \geq 20$.
This marked contrast with the relative encoding results quantitatively demonstrates the ability of our proposed encoding schemes to mitigate the curse of dimensionality.

\begin{table}[h]
 \tiny
  \caption{PROTES Naive Encoding without physically-constrained initialization. $N_{a}$: number of atoms; PC: PROTES calls (energy evaluations); Local calls total (LCT): total local optimization evaluations; Local calls last (LCL): last local optimization evaluations; RE: relative error; $SR_t$: success rate (\%).}
  \label{table:protes_naive_encoding}
  \centering
  \begin{tabular}{|c|r|r|r|r|r|}
    \hline
    $N_{a}$ & PC & LCT & LCL & RE & $SR_t$ \\
    \hline
5 & 1.27E+02 & 9.61E+03 & 3.23E+02 & 2.32E-12 & 100 \\
6 & 2.39E+03 & 1.16E+05 & 3.05E+02 & 2.05E-12 & 100 \\
7 & 1.39E+03 & 2.62E+04 & 3.13E+02 & 7.14E-11 & 100 \\
8 & 2.68E+02 & 9.00E+03 & 2.54E+02 & 2.18E-10 & 100 \\
9 & 8.66E+03 & 9.25E+03 & 2.39E+02 & 2.60E-11 & 100 \\
10 & 3.27E+04 & 1.00E+04 & 1.72E+02 & 2.81E-02 & 42 \\
11 & 2.75E+04 & 7.36E+03 & 1.68E+02 & 7.82E-03 & 78 \\
12 & 4.25E+04 & 3.18E+03 & 1.34E+02 & 9.55E-11 & 100 \\
13 & 3.85E+04 & 1.91E+03 & 2.56E+02 & 4.47E-02 & 50 \\
15 & 1.24E+05 & 7.64E+03 & 3.79E+02 & 1.28E-10 & 100 \\
16 & 1.80E+04 & 4.55E+03 & 3.96E+02 & 8.20E-11 & 100 \\
20 & 3.03E+05 & 6.69E+04 & 5.71E+02 & 3.06E-10 & 100 \\
    \hline
  \end{tabular}
\end{table}

\subsection{Complete Encoding Strategy Comparison}

Table~\ref{tab:encoding_comparison} synthesizes data across all encoding schemes and initialization strategies, providing a comprehensive guide for selecting the optimal configuration based on cluster size and performance priorities.
The results clearly underscore the superiority of relative encodings over direct encoding, as well as the systematic transition from agnostic to physically-constrained initialization as cluster size increases.

\begin{table}[H]
    \centering
   \tiny
    \caption{Comparison of optimization performance across different encoding schemes and initialization strategies for Lennard-Jones clusters. Values represent the total number of energy evaluations (including gradient calculations) required to locate the global minimum. Best results for each cluster size are highlighted in bold. Encoding schemes: SR = simpler relative encoding with constant distances; CR = relative encoding with spherical constraints; DE = direct Cartesian (naive) encoding. Initialization strategies: Agn = agnostic, uniform initialization; PhC = physically-constrained initialization. For direct Cartesian encoding, physically-constrained initialization is not applicable.}
    \label{tab:encoding_comparison}
    \begin{tabular}{|c|r|r|r|r|r|}
    \hline
    $N_a$ & \multicolumn{2}{c|}{SR} & \multicolumn{2}{c|}{CR} & \multicolumn{1}{c|}{DE} \\
    \cline{2-6}
    & Agn & PhC & Agn & PhC & Agn \\
    \hline
    5 & 1.39E+03 & \textbf{6.93E+02} & 3.29E+03 & 6.87E+02 & 9.74E+03 \\
    6 & 4.91E+03 & 5.68E+03 & \textbf{2.25E+03} & 2.34E+03 & 1.18E+05 \\
    7 & \textbf{5.48E+02} & 2.68E+03 & 1.99E+03 & 3.66E+03 & 2.75E+04 \\
    8 & \textbf{6.21E+02} & 1.45E+03 & 1.11E+03 & 2.46E+03 & 9.27E+03 \\
    9 & \textbf{9.66E+02} & 3.56E+03 & 1.28E+03 & 1.01E+04 & 1.79E+04 \\
    10 & 4.57E+03 & \textbf{3.11E+03} & 5.46E+03 & 4.32E+03 & - \\
    11 & \textbf{4.62E+03} & 1.15E+04 & 6.89E+03 & 7.02E+04 & - \\
    12 & 7.94E+03 & 4.55E+03 & - & \textbf{3.41E+03} & 4.57E+04 \\
    13 & 3.20E+03 & 5.54E+03 & \textbf{1.69E+03} & 7.42E+05 & - \\
    14 & 4.23E+03 & 4.01E+03 & \textbf{2.60E+03} & 3.85E+03 & - \\
    15 & 2.20E+03 & 3.99E+03 & \textbf{2.15E+03} & 2.00E+04 & 1.32E+05 \\
    16 & 8.91E+03 & 9.81E+03 & \textbf{6.51E+03} & 1.35E+04 & 2.26E+04 \\
    17 & 2.03E+04 & 3.37E+04 & \textbf{1.39E+04} & 1.87E+04 & - \\
    18 & 2.50E+05 & \textbf{1.57E+05} & - & 2.69E+05 & - \\
    19 & \textbf{1.21E+04} & 1.45E+04 & - & 1.66E+04 & - \\
    20 & 1.87E+04 & \textbf{1.32E+04} & - & 2.93E+04 & 3.70E+05 \\
    21 & 4.05E+04 & 1.18E+05 & 7.46E+04 & \textbf{2.82E+04} & - \\
    22 & 1.03E+05 & \textbf{6.72E+04} & - & 1.04E+05 & - \\
    23 & 1.15E+05 & - & - & \textbf{8.76E+04} & - \\
    24 & \textbf{8.24E+04} & 2.04E+05 & - & 7.95E+05 & - \\
    25 & 4.13E+04 & \textbf{1.18E+04} & - & 1.25E+04 & - \\
    26 & - & \textbf{2.34E+05} & - & 1.04E+06 & - \\
    27 & 2.92E+05 & \textbf{2.44E+05} & 3.30E+05 & 3.59E+05 & - \\
    28 & 1.48E+05 & \textbf{7.76E+04} & - & 2.98E+05 & - \\
    \hline
    \end{tabular}
\end{table}

\end{document}